\documentclass[reqno,11pt]{amsart}
\usepackage{amsmath,amsxtra,amsbsy,enumerate, latexsym, amsfonts, amssymb, amsthm, amscd, stmaryrd}
\usepackage{hyperref}
\usepackage{graphics,epsf,psfrag}
\usepackage{bbm, dsfont}
\usepackage{mathtools}
\usepackage{xcolor}
\usepackage{float}
\usepackage{tikz}
\usepackage{colonequals} 
\usepackage{microtype}

\numberwithin{equation}{section}
\allowdisplaybreaks
\newcommand\numberthis{\addtocounter{equation}{1}\tag{\theequation}}

\newcommand{\lint}{\llbracket}
\newcommand{\rint}{\rrbracket}

\newtheorem{theorem}{Theorem}[section]
\newtheorem{lemma}[theorem]{Lemma}
\newtheorem{proposition}[theorem]{Proposition}

\newtheorem{rem}[theorem]{Remark}

\newcommand{\Tau}{\mathcal{T}}
\newcommand{\TV}{\mathrm{TV}}
\newcommand{\T}{T_{\mathrm{mix}}^{L, \lambda}} 
\newcommand{\Tmaximal}{T_{\mathrm{mix}}^{L, \wedge}}
\newcommand{\Tminimal}{T_{\mathrm{mix}}^{L, \vee}} 
\newcommand{\Textremal}{\breve{T}_{\mathrm{mix}}^L}
\newcommand{\muL}{\mu_L^{ \lambda}} 
\newcommand{\barsin}{\overline{\sin}}
\newcommand{\barcos}{\overline{\cos}}
\newcommand{\barpsi}{\overline{\Psi}}
\newcommand{\barphi}{\overline{\Phi}}
\newcommand{\bargamma}{\gamma}
\newcommand{\bwedge}{\overline{\wedge}}
\newcommand{\barh}{\overline{H}}
\newcommand{\Var}{\mathrm{Var}}

\newcommand{\RN}[1]{%
  \textup{\uppercase\expandafter{\romannumeral#1}}%
}

\newcommand{\gap}{\mathrm{gap}}
\makeatletter
\def\captionfont@{\footnotesize}
\def\captionheadfont@{\scshape}

\long\def\@makecaption#1#2{%
  \vspace{2mm}
  \setbox\@tempboxa\vbox{\color@setgroup
    \advance\hsize-6pc\noindent
    \captionfont@\captionheadfont@#1\@xp\@ifnotempty\@xp
        {\@cdr#2\@nil}{.\captionfont@\upshape\enspace#2}%
    \unskip\kern-6pc\par
    \global\setbox\@ne\lastbox\color@endgroup}%
  \ifhbox\@ne 
    \setbox\@ne\hbox{\unhbox\@ne\unskip\unskip\unpenalty\unkern}%
  \fi
  \ifdim\wd\@tempboxa=\z@ 
    \setbox\@ne\hbox to\columnwidth{\hss\kern-6pc\box\@ne\hss}%
  \else 
    \setbox\@ne\vbox{\unvbox\@tempboxa\parskip\z@skip
        \noindent\unhbox\@ne\advance\hsize-6pc\par}%
\fi
  \ifnum\@tempcnta<64 
    \addvspace\abovecaptionskip
    \moveright 3pc\box\@ne
  \else 
    \moveright 3pc\box\@ne
    \nobreak
    \vskip\belowcaptionskip
  \fi
\relax
}
\makeatother

\def\writefig#1 #2 #3 {\rlap{\kern #1 truecm
\raise #2 truecm \hbox{#3}}}

\begin{document}
\title[Cutoff for Polymer pinning]{Cutoff for polymer pinning dynamics in the repulsive phase}

\author[Shangjie Yang]{Shangjie Yang}
 \address{Shangjie Yang \hfill\break
IMPA\\
Estrada Dona Castorina, 110\\ Rio de Janeiro 22460-320 \\ Brazil.}
\email{yashjie@impa.br}

\date{\small\today}

\begin{abstract}
  We consider the Glauber dynamics for model of polymer interacting with a substrate or wall. The state space is the set of one-dimensional nearest-neighbor paths on $\mathbb{Z}$ with nonnegative integer coordinates, starting at $0$ and coming back to $0$ after $L$ ($L\in 2\mathbb{N}$) steps and the Gibbs weight of a path $\xi=(\xi_x)^{L}_{x=0}$ is given by $\lambda^{\mathcal{N}(\xi)}$, where $\lambda \geq 0$ is a parameter which models the intensity of the interaction with the substrate and $\mathcal{N}(\xi)$ is the number of zeros in $\xi$. The dynamics we consider proceeds by updating $\xi_x$ with rate one for each $x=1,\dots, L-1$, in a heat-bath fashion. This model was introduced in  \cite{caputo2008approach} with the aim of studying the relaxation to equilibrium of the system.\\
 We present new results concerning the total variation mixing time for this dynamics when $\lambda< 2$,
 which corresponds to the phase where the effects of the wall's entropic repulsion dominates.
  For $\lambda \in [0, 1]$,
 we prove that  the total variation distance to equilibrium drops abruptly from $1$ to $0$ at time  $(L^2 \log L)(1+o(1))/\pi^2$.
 For $\lambda \in (1,2)$, we prove that the system also exhibit cutoff at time $(L^2 \log L)(1+o(1))/\pi^2$ when considering mixing time from ``extremal conditions'' (that is, either the highest or lowest initial configuration for the natural order on the set of paths).
 Our results improves both previously proved upper and lower bounds in \cite{caputo2008approach}.
\end{abstract}


\maketitle

\section{Introduction} 

\subsection{The random walk pinning model}
Consider the set of all one-dimensional nearest-neighbor paths on $\mathbb{Z}$ with nonnegative integer coordinates, starting at $0$ and coming back to $0$ after $L$ steps, \textit{i.e.} 
\begin{equation*}
\Omega_L\colonequals\Big\{\xi \in \mathbb{Z}^{L+1}  : \xi_0=\xi_L=0; |\xi_{x+1}-\xi_x|=1, \forall x \in \lint 0, L-1 \rint; \xi_x\geq 0,  \forall x\in \lint 0, L \rint  \Big\}, 
\end{equation*}
where $L\in 2\mathbb{N}$, and  $\lint a, b \rint \colonequals \mathbb{Z} \cap [a, b]$ for all $ a, b \in \mathbb{R}$ with $a<b$.  We study the polymer pinning model. This model is obtained by assigning to each path $\xi \in \Omega_L$ a weight 
$\lambda^{\mathcal{N}(\xi)}$, in which $\lambda \geq 0$ is the pinning parameter and
\begin{equation}
 \mathcal{N}(\xi)\colonequals\sum_{x=1}^{L-1}\mathbbm{1}_{\{\xi_x=0 \}}  \label{number of contact points}
\end{equation}
is the number of contact points with the $x$-axis. By convention,  $0^0 \colonequals 1$ and $0^n \colonequals 0$ for any positive integer $n\geq 1$.
Normalizing the weights, we obtain a Gibbs probability measure $\muL$ on $\Omega_L$, defined by
\begin{equation}
\muL(\xi)\colonequals\frac{\lambda^{\mathcal{N}(\xi)}}{Z_L( \lambda)} \label{the stationary measure}
\end{equation}
 where $\xi\in \Omega_L$ and
\begin{equation}
Z_L( \lambda)\colonequals\sum_{\xi' \in \Omega_L} \lambda^{\mathcal{N}(\xi')}.  \label{partition function}
\end{equation}
The graph of $\xi$ represents the spatial conformation of the polymer and $\lambda$ models the energetic interaction with an impenetrable substrate which fills the lower half plane ($\lambda<1$ corresponding to a repulsive interaction, $\lambda>1$ to an attractive one).
Since $\xi_x\geq 0$ for any $\xi\in \Omega_L$ and any $x\in \lint 0, L \rint$, we say that the polymers interact with an impenetrable substrate.
When there is no confusion, we drop the indices $\lambda$ and $L$ in $\muL$.

The random walk pinning model was  introduced in the seminal paper \cite{fisher1984walks} several decades ago,
and its various derivative models have been studied since. We refer to   \cite{GiacominPolymerbk, GiacominPolymerLNM} for recent reviews, and mention  \cite[Chapter 2]{GiacominPolymerbk} and references therein for more details. This model displays a transition from a delocalized phase to a localized phase  (see \cite[Section 1]{caputo2008approach}):  (a) if  $0\leq\lambda < 2$, the expected number of contact $\mu^{\lambda}_L(\mathcal N(\xi))$ is uniformly bounded in $L$ and the height of the middle point $\xi_{L/2}$ typically of order $\sqrt{L}$; (b) if $\lambda>2$, typical paths have a number of contact which is of order $L$ 
and distribution the heigh of the middle point $\xi_{L/2}$ is (exponentially) tight in $L$. These two phases are referred to as the delocalized and localized phase respectively, at the critical point $\lambda=2$ the system displays an intermediate behavior.

\medskip

A dynamical version of this model  was introduced more recently by Caputo \textit{et al.} in \cite{caputo2008approach}.  
     The corner-flip Glauber dynamics is a continuous-time reversible Markov chain on $\Omega_L$ with $\muL$ as the unique invariant probability measure, whose transitions are given by the updates of local coordinates. We refer to Figure \ref{fig:jumprates} for a graphical description of the jump rates for the system.     
     The dynamics is studied to understand how the system relaxes to equilibrium. Caputo \textit{et al.} in  \cite[Theorems 3.1 and 3.2]{caputo2008approach} proved that for $\lambda \in [0, 2)$, the mixing time of the dynamics in $\Omega_L$ is of order $L^2\log L$,  with non-matching constant prefactors for the upper and lower bounds.
     
     \medskip
 
  The goal of this paper is  to improve both the upper and lower bounds proved in \cite{caputo2008approach} and to show that the mixing time of the system is exactly $(1+o(1))(L^2 \log L)/\pi^2$ for $\lambda \in [0, 2)$. We prove the result for the worse initial condition mixing time when $\lambda\in [0,1]$. When $\lambda\in (1,2)$, our result is valid only for the mixing time starting from either the lowest or highest initial condition but we believe that this is only a technical restriction.

 \subsection{The dynamics} \label{Preliminaries}
\begin{figure}[h]
\centering
  \begin{tikzpicture}[scale=.4,font=\tiny]
   \draw (25,4) -- (25,-1) -- (52,-1);
    \draw[color=blue] (25,-1)--(26,0) -- (27,-1) -- (28,0) -- (29,1) -- (30,0) -- (31,-1) -- (32,0) -- (33,1) -- (34,0) -- (35,1) -- (36,2) -- (37,1) -- (38,0) -- (39,1) -- (40,0) -- (41,1) -- (42,0) -- (43,-1) -- (44,0) -- (45,-1) -- (46,0) -- (47,1) -- (48,2) -- (49,1) -- (50,0)--(51, -1);
    \foreach \x in {25,...,51} {\draw (\x,-1.3) -- (\x,-1);}
    \draw[fill] (25,-1) circle [radius=0.1];
    \draw[fill] (51,-1) circle [radius=0.1];   
    \node[below] at (25,-1.3) {$0$};
    \node[below] at (51,-1.3) {$L$};
     \draw[dashed] (33,1) -- (34,2) -- (35,1);
     \node[below] at (34.5, 2.8) {$\tfrac{1}{2}$};
      \draw (34, 0.3) edge[out=70, in=290, ->] (34, 1.8);
    \draw[dashed] (30,0) -- (31,1) -- (32,0);
    \draw (31,-.8) edge[out=70, in=290, ->] (31,.8);
     \node[below] at (31.7, 2) {$\tfrac{1}{1+\lambda}$};
    \draw[dashed, red](25,-1)--(26, -2)--(27,-1);
     \draw[red] (26,-0.2) edge[out=290, in=70, ->] (26,-1.8);
     \node[below, red] at (26.3,-0.2){$\times$};
 \draw[dashed, red](43,-1)--(44, -2)--(45,-1);  
   \draw[red] (44,-0.2) edge[out=290, in=70, ->] (44,-1.8);
     \node[below, red] at (44.3,-0.2){$\times$};
     
    \draw[dashed] (40,0) -- (41,-1) -- (42,0);
    \draw (41,.8) edge[out=290, in=70, ->] (41,-.8);
    \node[above] at (42,.4) {$\tfrac{\lambda}{1 + \lambda}$};
 \draw[dashed] (47,1) -- (48,0) -- (49,1);
    \draw (48,1.8) edge[out=290, in=70, ->] (48,0.2);
    \node[below] at (48.5, 2.7) {$\tfrac{1}{2}$};
    \node[blue,above] at (38, 1.8){$\xi$};
    \draw[thick,->] (25,-1) -- (25,4) node[anchor=north west]{y};
    \draw[thick,->] (25,-1) -- (52,-1) node[anchor=north west]{x};
  \end{tikzpicture}
  \caption{\label{fig:jumprates}  A graphical representation of the jump rates for the system pinned at $(0, 0)$ and $(L, 0)$.  
  A transition of the dynamics corresponds to flipping a corner, whose rate depends on how it changes the number of contact points with the $x$-axis. The rates are chosen in a manner such that the dynamics is reversible with respect to $\muL$. The two red dashed corners are not available and labeled with $\times$, because of the nonnegative restriction of the state space $\Omega_L$. Note that not all the possible transitions are shown in the figure. }
\end{figure}
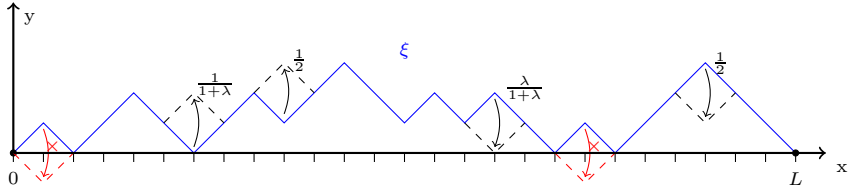

For $\xi \in \Omega_L$ and $x \in \lint 1, L-1\rint$,  we define $\xi^x\in \Omega_L$ by
\begin{equation}\label{flippingcorner}
\xi^x_y\colonequals\begin{cases*}
\xi_y &if $y\neq x$,\\
(\xi_{x-1}+\xi_{x+1})-\xi_x & if $y=x \mbox{ and }  \xi_{x-1}=\xi_{x+1}\geq 1 \mbox{ or } \xi_{x-1}\neq \xi_{x+1}$,\\
\xi_x  & if  $y=x \mbox{ and } \xi_{x-1}=\xi_{x+1}=0$.
\end{cases*}
\end{equation}
When $\xi_{x-1}=\xi_{x+1}$, $\xi$ displays a local extremum at $x$ and we obtain $\xi^x$ by flipping the corner of $\xi$ at the coordinate $x$,  provided that the path obtained by flipping the corner is in $\Omega_L$ (this excludes corner-flipping when  $\xi_{x-1}=\xi_{x+1}=0$.). See Figure \ref {fig:jumprates} for a graphical representation.
Given the probability measure $\muL$ defined in (\ref{the stationary measure}), we construct a continuous-time Markov chain whose generator $\mathcal{L}$ is given by its action on the functions $\mathbb{R}^{\Omega_L}$.  It can be written explicitly as
\begin{equation}
(\mathcal{L}f)(\xi) \colonequals\sum_{x=1}^{L-1} R_x(\xi) \big[f(\xi^x) -f(\xi) \big],  \label{the generator of the dynamics}
\end{equation}
where $f \colon \Omega_L \to\mathbb{R}$ is a function, and
\begin{equation*}
    R_x (\xi): =   \begin{cases*}
      \frac{1}{2}  & if $\xi_{x-1}=\xi_{x+1}>1$, \\
       \frac{\lambda}{1+\lambda}        &  if $(\xi_{x-1}, \xi_x, \xi_{x+1})=(1, 2, 1) $,\\
       \frac{1}{1+\lambda}        &  if $(\xi_{x-1}, \xi_x, \xi_{x+1})=(1, 0, 1) $,\\
       0 & if $\xi_{x-1}\neq \xi_{x+1}$ or $\xi_{x-1}=\xi_{x+1}=0$.
    \end{cases*}
\end{equation*}

 \noindent Equivalently, we can rewrite the generator as
\begin{equation}
(\mathcal{L}f)(\xi)=\sum_{x=1}^{L-1} \Big[Q_x(f)(\xi)-f(\xi)\Big], \label{the generator of the dynamics2}
\end{equation}
and 
\begin{equation*}
Q_x(f)(\xi)\colonequals\muL \Big(f \big\vert  (\xi_y)_ {y \neq x}\Big).
\end{equation*}
Let $(\sigma_t^{\xi, \lambda})_{t\geq 0}$ be the trajectory of the Markov chain with initial condition $\sigma_0^{\xi}=\xi$ and parameter $\lambda$, and let $P_t^{\xi, \lambda}$ be the law of distribution of the time marginal $\sigma_t^{\xi, \lambda}$. 
Since $\muL(\xi) R_x(\xi)=\muL(\xi^x) R_x(\xi^x)$, 
the continuous-time chain is reversible 
 with respect to the probability measure  $\muL$. This chain is called the Glauber dynamics. Because the Markov chain is irreducible,  by  \cite[ Theorem 3.5.2]{Norris} we know that for all $\xi \in \Omega_L$, $P_t^{\xi, \lambda}$ converges to $\muL$ in discrete topology as $t$ tends to infinity. 
 We ask a quantitative question: how long does it take for $P_t^{\xi, \lambda}$ to converge to $\muL$, especially for the worst initial starting path $\xi \in \Omega_L$?

Let us state the aforementioned question in a mathematical framework. If $\alpha$ and $\beta$ are two probability measures on the space $(\Omega_L, 2^{\Omega_L})$, the total variation distance between $\alpha$ and $\beta$ is 
\begin{equation}
\Vert\alpha-\beta\Vert_{\TV}\colonequals\frac{1}{2}\sum_{\xi \in \Omega_L} \vert \alpha(\xi) -\beta(\xi) \vert=\sup_{\mathcal{A}\subset \Omega_L} \Big(\alpha(\mathcal{A})-\beta(\mathcal{A})\Big). \label{definition of total variation distance}
\end{equation}
We define the distance to equilibrium at time $t$ by
\begin{equation}
d^{L, \lambda}(t)\colonequals\max_{\xi \in \Omega_L} \Vert P^{\xi, \lambda}_t -\muL \Vert_{\TV}.
\end{equation}  
For any given $\epsilon\in (0, 1)$, let  the $\epsilon$-mixing-time  be
\begin{equation}
\T(\epsilon)\colonequals\inf \big\{ t\geq 0 : d^{L, \lambda}(t)\leq \epsilon  \big\}.
\end{equation}
We say that this sequence of Markov chains has a cutoff, if for all $\epsilon\in (0, 1)$,
\begin{equation}
\lim_{L\rightarrow \infty} \frac{\T(\epsilon)}{\T(1-\epsilon)}=1.
\end{equation}
The cutoff phenomenon is surveyed in the seminal paper \cite{diaconis1996cutoff},  and we refer to \cite[Chapter 18]{LPWMCMT} for more information.
In \cite[Theorems 3.1 and 3.2]{caputo2008approach}, for all $\lambda\in [0, 2)$,  Caputo \textit{et al.} proved that for all $\delta>0$ and all $\epsilon \in (0, 1)$, if $L$ is sufficiently large, we have 
\begin{equation}
\frac{1-\delta}{2\pi^2}L^2\log L \leq  \T(\epsilon) \leq \frac{6+\delta}{\pi^2}L^2\log L.
\end{equation}

\subsection{Main results}

In this paper, we find that the mixing time is $\big(1+o(1)\big)(L^2 \log L)/\pi^{2}$ for all $\lambda \in [0, 1]$, improving both the lower and upper bounds in \cite{caputo2008approach}. That is the following theorem.
\begin{theorem}\label{the cutoff phenomenon in the diffusive regime}
 For all $\epsilon \in (0, 1)$ and all $\lambda\in [0, 1]$,  we have
  \begin{equation}
  \lim_{L\rightarrow \infty} \frac{\pi^2 \T(\epsilon)}{L^2 \log L}=1.
  \end{equation}
\end{theorem}
Therefore, there is a cutoff phenomenon in the Glauber dynamics for $\lambda \in [0, 1]$. The reason why we include the result for $\lambda=0$ is the need for the mixing time about the dynamics when $\lambda \in (1, 2)$. 

\begin{rem}\label{remark}
Theorem \ref{the cutoff phenomenon in the diffusive regime} about $\lambda=0$ is the same as the case $\lambda=1$ by the following identification.
Let 
\begin{equation}
\Omega_L^{+}\colonequals \Big\{\xi \in \Omega_L: \mathcal{N}(\xi)=0 \Big\}
\end{equation}
where $\mathcal{N}(\xi)$ is defined in (\ref{number of contact points}),
and identify $\Omega_L^{+}$ with $\Omega_{L-2}$  by lifting up the $x$-axis up by distance one in $\Omega_L$. Precisely,
the identification is as follows: $\xi=(\xi_x)_{x\in \lint 0, L \rint} \in \Omega_L^{+}$ is identified with $\varsigma=(\varsigma_x)_{x\in \lint 0, L-2 \rint}\in \Omega_{L-2}$,  if $\varsigma_x=\xi_{x+1}-1$ for all $x \in \lint 0, L-2\rint$.
We can see:
\begin{itemize}
\item[(a)] $\mu_L^0$ is the same as the probability measure $\mu_{L-2}^{1}$;
\item[(b)] the dynamics $(\sigma_t^{\xi, 0})_{t\geq 0}$\textemdash living in the space $\Omega^+_L$\textemdash is the same as the dynamics $(\sigma_t^{\varsigma,1})_{t\geq 0}$ living in the space $\Omega_{L-2}$, where $\xi \in \Omega_L^+$ is identified with $\varsigma\in \Omega_{L-2}$.
 \end{itemize}
\end{rem}

 Therefore, we only need to prove Theorem \ref{the cutoff phenomenon in the diffusive regime} for $\lambda \in (0, 1]$.
In addition, we have a partial result for $\lambda \in (1, 2)$. Let us state the framework.
We introduce a natural partial order $``\leq"$ on $\Omega_L$ as follows
\begin{equation*}
\Big(\xi \leq \xi'\Big) \Leftrightarrow \Big(\forall  x\in \lint 0, L \rint, \xi_x \leq \xi'_x\Big).
\end{equation*}
In other words, if $\xi\leq \xi'$, the path $\xi$ lies below the path $\xi'$. Then the maximal path $\wedge$ and the minimal path $\vee$ are respectively given by
\begin{align*}
  \wedge_x : &= \min \big( x, -x+L \big), \quad \forall x \in \lint 0, L \rint; \\   
      \vee_x : &= x- 2 \lfloor x/2 \rfloor, \quad \forall x \in \lint 0, L \rint.   
\end{align*}
where $\lfloor x/2 \rfloor \colonequals \sup \big\{n\in \mathbb{Z}: n\leq x/2 \big\}$.
Define
\begin{align*}
\Tmaximal(\epsilon)&\colonequals\inf \Big\{ t\geq 0 :  \Vert P_t^{\wedge,\lambda} -\muL\Vert_{\TV} \leq \epsilon \Big\},\\
\Tminimal(\epsilon)&\colonequals\inf \Big\{ t\geq 0 :  \Vert P_t^{\vee,\lambda} -\muL\Vert_{\TV} \leq \epsilon \Big\},
\end{align*}
and 
\begin{equation}
\Textremal(\epsilon)  \colonequals \max \Big( \Tmaximal, \Tminimal \Big).
\end{equation}
For $ \lambda \in (1, 2)$, applying Peres-Winkler censoring inequality in \cite[Theorem 1.1]{peres2013can}, 
 we discover that the mixing time is also  $(1+o(1))(L^2 \log L)/\pi^2$ for the dynamics starting with the two extremal paths. That is the following theorem.
\begin{theorem}\label{the cutoff for the extremal paths}
For all $\epsilon \in (0, 1)$ and $\lambda \in (1, 2)$, we have
\begin{equation}
 \lim_{L\rightarrow \infty} \frac{\pi^2 \Textremal(\epsilon)}{L^2 \log L}=1.
\end{equation}
\end{theorem}

\subsection{Other values of $\lambda$}
Our analysis excludes the case $\lambda> 2$, let us just mention that the 
 convergence to equilibrium follows a different pattern in this case.
 While the relaxation time and and mixing time are of order 
 $L^{2}$ and $L^2\log L$ in the repulsive phase $\lambda <2$, it is believed that 
 they become of order $L$ and $L^2$ respectively in the attractive phase.
Rigorous lower bound has been proved in  \cite[Theorem 3.2]{caputo2008approach}, but matching order upper bound has only been shown when $\lambda=\infty$ (\cite[Proposition 5.6]{caputo2008approach} for the mixing time).
Furthermore in \cite[Theorem 2.7]{lacoin2014scaling}, it is shown that in this last case mixing time is equal to $L^2/4(1+o(1))$. When
$\lambda\in (2,\infty)$,  the conjecture in 
\cite[Section 2.7]{lacoin2014scaling} seems to indicate that the mixing time should be of order $C(\lambda) L^2(1+o(1))$ for some explicit $C(\lambda)$.

\medskip

At the critical value $\lambda=2$, we believe that the mixing time continues to be $\frac{L^2}{\pi^2}(\log L)(1+o(1))$ but our techniques do not allow to treat this case.

\subsection{Organization of the paper}
Section \ref{technical-preminary} introduces a grand coupling for the dynamics corresponding to different $\xi$ and $\lambda$, and some useful reclaimed results.

Section \ref{lower bounf for the mixing time} is dedicated to the lower bound on the mixing time for $\lambda \in (0, 2)$.

Section \ref{upper bound for the mixing time in the first diffusive regime} supplies the upper bound on the mixing time for $\lambda \in (0, 1]$.

Section \ref{upper bound for the mixing time of the dynamics starting with extremal paths} is about the upper bound on the mixing time for the dynamics starting with the two extremal paths when $\lambda \in (1, 2)$, applying censoring inequality.

\subsection{Notation}   
          We use $``\colonequals"$ to define a new  quantity on the left-hand side, and use $``\equalscolon"$ in some cases when the
quantity is defined on the right-hand side. 
 
    We let $(C_n(\lambda))_{n\in \mathbb{N}}$ and $(c_n(\lambda))_{n \in \mathbb{N}}$ be some positive constants, which are only dependent on $\lambda$. Additionally, we let  $(c_n)_{n\in \mathbb{N}}$ and $(C_n)_{n \in \mathbb{N}}$ be some positive and universal constants.

\subsection{Acknowledgments}I'm very grateful to my Ph.D supervisor Hubert Lacoin for suggesting this problem, inspiring discussions and comments on the manuscript. In addition, I also thank Anna Ben-Hamou and Daniela Cuesta for reading the manuscript carefully and their suggestions
for improving the manuscript.
 
     \medskip

\section{Technical preliminaries}\label{technical-preminary}
To use the monotonicity of the Glauber dynamics, we provide a graphical construction of the Markov chain such that all dynamics, \textit{i.e.} $\big\{ (\sigma_t^{\xi, \lambda})_{t\geq 0}: \forall \xi \in \Omega_L, \forall \lambda \in [0, \infty) \big\}$, live in one common probability space. 
This  construction appears in \cite[Section 8.1]{lacoin2016mixing}, which provides more independent flippable corners by comparing with the coupling in \cite[Subsection 2.2.1]{caputo2008approach}.
\subsection{A graphical construction.} \label{graphical construction} We set the exponential clocks and independent ``coins" in the centers of the squares formed by all the possible corners and their counterparts.
Let
\begin{equation}
\Theta\colonequals\Big\{(x, z) :  x\in \lint 2, L-2 \rint,  z\in \lint 1, L/2-1-\vert x- L/2 \vert \rint; x+z\in 2\mathbb{N}+1 \Big\}, \label{the set of spins for censoring}
\end{equation}
and let $\Tau^{ \uparrow}$ and $\Tau^{\downarrow}$ be two independent rate-one exponential clock processes indexed by $\Theta$. 
That is to say,  for every $(x, z)\in \Theta$ and $n\geq 0$, we have $\Tau_{(x,z)}^{\uparrow}(0)=0$, and 
 \begin{equation*}
 \Big( \Tau_{(x, z)}^{\uparrow}(n)-\Tau_{(x, z)}^{\uparrow}(n-1) \Big)_{n\geq 1}
 \end{equation*}
 is a field of i.i.d. exponential random variables with mean $1$. Similarly, this holds for $\Tau_{(x, z)}^{\downarrow}$. Moreover, let $\mathcal{U}^{\uparrow}=\Big(U^{\uparrow}_{(x, z)}(n)\Big)_{(x, z)\in \Theta, n\geq 1}$ and $\mathcal{U}^{\downarrow}=\Big(U^{\uparrow}_{(x, z)}(n)\Big)_{(x, z)\in \Theta, n\geq 1}$ be two independent fields of i.i.d.  random variables uniformly distributed in $[0, 1]$, which are independent of $\Tau^{\uparrow}$ and $\Tau^{\downarrow}$. 
Given $\Tau^{\uparrow}$, $\Tau^{\downarrow}$, $\mathcal{U}^{\uparrow}$ and $\mathcal{U}^{\downarrow}$,  we construct, in a deterministic way,  $(\sigma_t^{\xi, \lambda})_{t \geq 0}$  the trajectory of the Markov chain with parameter $\lambda$ and starting with $\xi \in \Omega_L$, \textit{i.e.} $\sigma_0^{\xi, \lambda}=\xi$.

When the clock process $\Tau^{\uparrow}_{(x, z)}$ rings at time $t=\Tau_{(x, z)}^{\uparrow}(n)$  for $n\geq 1$ and $\sigma_{t^-}^{\xi, \lambda}(x)=z-1$,  we update $\sigma_{t^-}^{\xi, \lambda}$ as follows:
\begin{itemize}
\item if $\sigma_{t^-}^{\xi,\lambda}(x-1)=\sigma_{t^-}^{\xi,\lambda}(x+1)=z\geq 2$ and $U_{(x, z)}(n)^{\uparrow}\leq \frac{1}{2}$, let $\sigma_{t}^{\xi, \lambda}(x)=z+1$ and the other coordinates remain  unchanged;

\item if $\sigma_{t^-}^{\xi,\lambda}(x-1)=\sigma_{t^-}^{\xi, \lambda}(x+1)=z=1$ and $U_{(x, z)}^{\uparrow}(n)\leq \frac{1}{1+\lambda}$, let $\sigma_{t}^{\xi, \lambda}(x)=2$ and the other coordinates remain  unchanged.

\end{itemize}
 If neither of these two aforementioned conditions is satisfied, we do nothing.

When the clock process $\Tau^{\downarrow}_{(x, z)}$ rings at time $t=\Tau_{(x,z)}^{\downarrow}(n)$ for $n\geq 1$ and $\sigma_{t^-}^{\xi, \lambda}(x)=z+1$, we update $\sigma_{t^-}^{\xi,\lambda}$ as follows:
\begin{itemize}
\item if $\sigma_{t^-}^{\xi,\lambda}(x-1)=\sigma_{t^-}^{\xi,\lambda}(x+1)=z\geq 2$ and $U_{(x, z)}^{\downarrow}(n)\leq \frac{1}{2}$, let $\sigma_{t}^{\xi,\lambda}(x)=z-1$ and the other coordinates remain  unchanged;

\item if $\sigma_{t^-}^{\xi,\lambda}(x-1)=\sigma_{t^-}^{\xi,\lambda}(x+1)=z-1=0$ and $U_{(x, z)}^{\downarrow}(n)\leq \frac{\lambda}{1+\lambda}$, let $\sigma_{t}^{\xi, \lambda}(x)=0$ and the other coordinates remain  unchanged.
\end{itemize}
 If neither of these two aforementioned conditions is satisfied, we do nothing.

Let $\mathbb{P}$ or $\mathbb{E}$ stand for the probability law corresponding to $\Tau^{\uparrow}$, $\Tau^{\downarrow}$, $\mathcal{U}^{\uparrow}$ and $\mathcal{U}^{\downarrow}$. 
Recall that  $\mu$ is the stationary probability measure for the dynamics.
The dynamics  $(\sigma_t^{\mu, \lambda})_{t\geq 0}$ is constructed by first taking the initial path $\xi$ sampling from $\mu$ at $t=0$ and then
using the graphical construction with parameter $\lambda$ for $t>0$. This sampling is independent of $\mathbb{P}$.
Define $P^{\mu, \lambda}_t(\cdot)\colonequals\mathbb{P}(\sigma^{\mu, \lambda}_t=\cdot)$, and likewise $P^{\mu, \lambda}_t(\mathcal{A})\colonequals\mathbb{P}[\sigma_t^{\mu, \lambda}\in \mathcal{A}]$ for $\mathcal{A}\subset \Omega_L$. 
 When it is clear in the context, we use the notations $(\sigma_t^{\mu})_{t\geq 0}$ and $P^{\mu}_t$, ignoring the parameter $\lambda$.

This graphical construction allows us to construct all the trajectories $(\sigma_t^{\xi, \lambda})_{t\geq 0}$ starting from all $\xi \in \Omega_L$ and all parameters $\lambda \in [0, \infty)$ simultaneously.
It  preserves the order, affirmed in the following proposition. The proof of this  proposition, which we omit, is almost identical to that of \cite[Proposition 3.1]{lacoin2016mixing}.

\begin{proposition}\label{preserving the monotonicity}
Let $\xi$ and $\xi'$ be two elements of $\Omega_L$ satisfying $\xi\leq \xi'$, and $0\leq \lambda\leq \lambda'$.  With the graphical construction above, we have
\begin{equation}
\begin{aligned}
&\mathbb{P}\Big[\forall t \in [0, \infty): \sigma_t^{\xi, \lambda}\leq \sigma_t^{\xi', \lambda}\Big]=1,\\
&\mathbb{P}\Big[\forall t \in [0, \infty):  \sigma_t^{\xi, \lambda'} \leq \sigma_t^{\xi, \lambda}\Big]=1. \label{monotonicity of the dynamics equations}
\end{aligned}
\end{equation}
\end{proposition}

\subsection{Useful reclaimed results.}
 We have the asymptotic information about $Z_L(\lambda)$, which is: 
\begin{theorem}[Theorem 2.1 in \cite{caputo2008approach}] \label{theorem for the asymptotic behavior of the partition function}
For every $\lambda \in [0, 2)$, we have
\begin{equation}
\lim_{L \to \infty} \frac{Z_L(\lambda)}{2^{L}L^{-3/2}}=C(\lambda),
\end{equation}
where $C(\lambda)>0$ is a constant, only dependent on $\lambda$.
\end{theorem}

Furthermore,
to understand the Glauber dynamics, it is important to understand how  the generator $\mathcal{L}$ acts on the paths in $\Omega_L$. Let us introduce the settings.
For a function $g\colon \lint 0, L\rint \to \mathbb{R}$, the discrete Laplace operator $\Delta$ is defined as follows:  for any $x\in \lint 1 , L-1 \rint$,
$$(\Delta g)_x\colonequals\frac{1}{2}\Big(g(x-1)+g(x+1)\Big)-g(x).$$
Besides, we define a function $f\colon \Omega_L \to \mathbb{R}$ to be $f(\xi)\colonequals\xi_x$, and let $\mathcal{L}\xi_x\colonequals (\mathcal{L}f)(\xi)$ for  $x\in \lint 1 , L-1 \rint$. Considering (\ref{the generator of the dynamics2}), we know that $\mathcal{L}\xi_x=\muL(\xi_x \vert \xi_{x-1}, \xi_{x+1})-\xi_x$, and a calculation yields the following identity which we recall as a lemma.
\begin{lemma}[Lemma 2.3 in \cite{caputo2008approach}]
\label{the action of operator on paths}    For every $\lambda>0$ and every $ x \in \lint 1, L-1 \rint$, we have
\begin{equation}
\mathcal{L}\xi_x=(\Delta \xi)_x+\mathbbm{1}_{\{ \xi_{x-1}=\xi_{x+1}=0 \}}-\bigg(\frac{\lambda-1}{\lambda+1}\bigg)\mathbbm{1}_{\{\xi_{x-1}=\xi_{x+1}=1\}}.
\end{equation}
\end{lemma}

\section{Lower bound on the mixing time for $\lambda\in (0, 2)$}
\label{lower bounf for the mixing time}
This section is devoted to providing a lower bound on the mixing time of the Glauber dynamics for $\lambda \in (0, 2)$, which is the following proposition.
\begin{proposition}\label{lower bound of the mixing time} 
For all $\lambda \in (0, 2)$ and all $\epsilon \in (0, 1)$, we have
\begin{equation}
\T(\epsilon)\geq \tfrac{1}{\pi^2}L^2\log L-C(\lambda, \epsilon)L^2\equalscolon t_{C(\lambda, \epsilon)}, \label{lower bound on the mixing time: inequality}
\end{equation}
where $C(\lambda, \epsilon)>0$ is a constant, only dependent on $\lambda$ and $\epsilon$. 
\end{proposition}

Before we start the proof of Proposition \ref{lower bound of the mixing time}, let us explain the idea. Note that
the function $\Phi(\xi)$, defined in (\ref{weighted area function introduced by Caputo}) below,  is almost the area enclosed by  the $x$-axis and the path $\xi \in \Omega_L$.
Intuitively, $\Phi(\wedge)$ is of order $L^2$, while at equilibrium $\Phi(\xi)$ is of order $L^{3/2}$.
 The second moment method in \cite[Theorem 3.2]{caputo2008approach} does not supply a sharp lower bound on the mixing time.
 We adapt the idea in \cite[Theorem 3.2]{caputo2008approach} to provide the lower bound in (\ref{lower bound on the mixing time: inequality}) by  proving the following.
 \begin{itemize}
 \item[(i)] While the expected equilibrium value $\mu(\Phi)$ is at most of order $L^{3/2}$, $\mathbb{E}[\Phi(\sigma_t^{\wedge})]$ is much bigger than $L^{3/2}$ for all $t\leq t_{C(\lambda, \epsilon)}$; 
 \item[(ii)] On the one hand $\Phi(\sigma_t^{\mu})$  is fairly close to its mean $\mu(\Phi)$ by Markov's inequality,  and on the other hand $\Phi(\sigma_t^{\wedge})$ is well concentrated around $\mathbb{E}[\Phi(\sigma_t^{\wedge})]$ by controlling its fluctuation through martingale approach.
 \end{itemize}
Subsection \ref{ingredientsforlowerbound} prepares all the ingredients for the first step of this strategy, and Subsection \ref{proofforlowerbound} is dedicated to the second step of of the strategy, giving the proof of Proposition \ref{lower bound of the mixing time}.

\subsection{ Ingredients for the lower bound of the mixing time. }\label{ingredientsforlowerbound}
Inspired by \cite[Equation (1)]{wilson2004mixing}, Caputo \textit{et al.} in \cite[Equation (2.39)]{caputo2008approach} defined
the weighted area function $\Phi\colon \Omega_L \to \mathbb{R}$ by
\begin{equation}
\Phi(\xi)\colonequals\sum_{x=1}^{L-1} \xi_x \barsin (x), \label{weighted area function introduced by Caputo}
\end{equation}
where $\barsin (x)\colonequals\sin(\frac{\pi x}{L})$ and $\xi \in \Omega_L$. As \cite[Equation (4.3)]{caputo2008approach},  we use Lemma \ref{the action of operator on paths} and  summation by part to obtain
\begin{equation}
(\mathcal{L}\Phi)(\xi)=\sum_{x=1}^{L-1} \barsin(x) \mathcal{L}\xi_x=-\kappa_L \Phi(\xi)+\Psi(\xi), \label{operator acting on weighted path}
\end{equation}
 where $\kappa_L\colonequals1-\cos (\frac{\pi}{L})$ and  
\begin{equation}
 \Psi(\xi)\colonequals\sum_{x=1}^{L-1}\barsin(x) \bigg[\mathbbm{1}_{\{ \xi_{x-1}=\xi_{x=1}=0  \}}-\bigg(\frac{\lambda-1}{\lambda+1}\bigg)\mathbbm{1}_{\{\xi_{x-1}=\xi_{x+1}=1 \}}\bigg]. \label{the effect of the wall}
\end{equation}
Since $\barsin(x)\geq 0$ for all $x \in \lint 0, L\rint$, we have
\begin{equation}
 |\Psi(\xi)| \leq \sum_{x=1}^{L-1}\barsin(x) \bigg[\mathbbm{1}_{\{ \xi_{x-1}=\xi_{x=1}=0  \}}+\bigg|\frac{\lambda-1}{\lambda+1}\bigg|\mathbbm{1}_{\{\xi_{x-1}=\xi_{x+1}=1 \}}\bigg] \equalscolon \barpsi(\xi). \label{absolute effect of the wall}
\end{equation}
Caputo \textit{et al.} gave an upper bound on $\muL(\Phi)$. In \cite[Equation (5.15)]{caputo2008approach}, they used coupling and monotonicity  to obtain that for every positive integer $k$,
\begin{equation*}
\sup_{\lambda\geq 0, L} \hspace{0.05cm} \sup_{x \in \lint 1, L-1 \rint } \muL \bigg(\frac{(\xi_x)^k}{L^{k/2}}\bigg)<\infty.
\end{equation*}
Consequently, using $k=1$ and $\barsin(x)\leq 1$,  we have
 \begin{equation}
 \muL(\Phi)\leq \sum_{x=1}^{L-1}\mu_L^{\lambda}(\xi_x)\leq c L^{3/2}, \label{upper bound for the weighted area}
 \end{equation}
  where $c>0$ does not depend on $\lambda$.
In addition, Caputo \textit{et al.} also gave a lower bound on $\mathbb{E}[\Phi(\sigma_t^{\wedge})]$, which we recall as a lemma below. 
\begin{lemma} [Equation (5.24) in \cite{caputo2008approach}] \label{lower bound of the area of the maximal path} For all $\lambda\in (0, 2)$,  all $t\geq 0$, all $L\geq 2$ and some constant $c(\lambda)>0$, we have
\begin{equation*}
\mathbb{E}[\Phi(\sigma_t^{\wedge})]\geq \Phi(\sigma_0^{\wedge}) e^{ -\kappa_L t}-c(\lambda)L^{3/2}. 
\end{equation*} 
\end{lemma}
In view of (\ref{absolute effect of the wall}), we need an upper bound on $\mathbb{P}\big[\sigma_t^{\wedge}(x-1)=\sigma_t^{\wedge}(x+1)\in \{0, 1\} \big]$ for $x \in \lint 1, L-1\rint$, which is the following lemma.

\begin{lemma}\label{lemma for the maximal dynamics hitting 0 or 1}
 For all $t\geq 0$, all $x\in \lint 1, L-1\rint$ and all $L\geq 2$, we have
\begin{equation}
\mathbb{P}\big[\sigma^{\wedge}_t(x-1)=\sigma_t^{\wedge}(x+1) \in \{0, 1 \}\big] \leq C_1(\lambda)\frac{L^{3/2}}{x^{3/2}(L-x)^{3/2}}. \label{upper bound for hitting 0 or 1 by monotonicity}
\end{equation}
\end{lemma}
\begin{proof}
Since $\sigma^{\wedge}_t\geq \sigma^{\mu}_t$ for all $t\geq 0$, 
we know that for all $x\in \lint 1, L-1\rint$,  
\begin{align*}
 \mathbb{P}\big[\sigma^{\wedge}_t(x-1)=\sigma_t^{\wedge}(x+1)\in \{ 0, 1\}\big]&\leq \mathbb{P}\big[\sigma^{\mu}_t(x-1)=\sigma_t^{\mu}(x+1)\in \{ 0, 1\}\big]\\
 &= \muL\big(\xi_{x-1}=\xi_{x+1}\in \{ 0, 1\}\big).   
\end{align*}
For all $\lambda\in (0, 2)$, all $x \in \lint 1, L-1 \rint \cap 2\mathbb{N}$ and all $L\geq 2$, applying Theorem \ref{theorem for the asymptotic behavior of the partition function},  we obtain 
\begin{equation}
\begin{aligned}
\muL(\xi_x=0)&=\lambda \frac{Z_x(\lambda) Z_{L-x}(\lambda)}{Z_L(\lambda)}
\leq C_2(\lambda)\frac{L^{3/2}}{x^{3/2}(L-x)^{3/2}}  \label{estimation for a hitting of the origin},
\end{aligned}
\end{equation}
and 
\begin{equation}
\muL(\xi_{x-1}=\xi_{x+1}=0)=\lambda^2\frac{Z_{x-1}(\lambda) Z_{L-x-1}(\lambda)}{Z_L(\lambda)}\leq C_2(\lambda)\frac{L^{3/2}}{x^{3/2}(L-x)^{3/2}}. \label{estimation of hitting the origin consecutively}
\end{equation}
With the same conditions about $\lambda$, $x$ and $L$ above,  as \cite[Equation (5.23)]{caputo2008approach} we have
\begin{equation}
\muL(\xi_{x-1}=\xi_{x+1}=1)=\frac{1+\lambda}{\lambda}\muL(\xi_x=0).
\label{estimation for a pair of point hitting of 1}
\end{equation}
Therefore, by (\ref{estimation for a hitting of the origin}), (\ref{estimation of hitting the origin consecutively}) and (\ref{estimation for a pair of point hitting of 1}), we obtain  (\ref{upper bound for hitting 0 or 1 by monotonicity}).
\end{proof}

\subsection{Proof of the lower bound on the mixing time.}\label{proofforlowerbound}
Let us  detail the second step of the aforementioned strategy.  To prove that $\Phi(\sigma_t^{\wedge})$ is well concentrated around its mean $\mathbb{E}[\Phi(\sigma_t^{\wedge})]$, we do the following.
\begin{itemize}
\item[(i)] For a fixed time $t_0$, we use the function $F(t,\xi)=\exp(\kappa_L(t-t_0))\Phi(\xi)$ to construct a Dynkin's  martingale $M$ (see \cite[Lemma 5.1 in Appendix 1]{LandimHydrodynamicsbk}).
\item[(ii)] To estimate the fluctuation of $F(t_0, \sigma_{t_0}^{\wedge})=\Phi(\sigma_{t_0}^{\wedge})$, we control the martingale bracket $\langle M. \rangle$ and the mean of $(\partial_t+\mathcal{L})F(t, \sigma_t^{\wedge})$, which comes from the construction of  Dynkin's martingale.
\end{itemize}
 While $\Phi(\sigma_t^{\mu})$ is at most of order $L^{3/2}$, $\Phi(\sigma_{t_0}^{\wedge})$ is much bigger than $L^{3/2}$ for all $t_0\leq t_{C(\lambda, \epsilon)}$. This property of $\Phi$ about $\sigma_t^{\mu}$ and $\sigma_{t_0}^{\wedge}$ can be used to provide a lower bound on the distance between $\mu$ and $P_{t_0}^{\wedge}$. 

\begin{proof}[Proof of Proposition \ref{lower bound of the mixing time}.]
We adapt the approach in  \cite[Proposition 5.3]{caputo2008approach}. For $C \in (0, \infty)$, define
\begin{equation}
\mathcal{A}_C \colonequals \big\{\xi \in \Omega_L: \Phi(\xi)\leq CL^{3/2} \big\}.\label{definition of the event AC}
\end{equation}
Using Markov's inequality and  (\ref{upper bound for the weighted area}), we obtain
\begin{equation}
1-\mu(\mathcal{A}_C)=\mu(\Phi > CL^{3/2})\leq \frac{\mu(\Phi)}{CL^{3/2}}\leq \frac{c}{C}, \label{lower bound for the well chosen event of the equilibrium states for the size of the windows}
\end{equation}
where the rightmost term is smaller than or equal to $\epsilon/2$ for $C\geq 2c/\epsilon$.
Our next step is to prove that for any given $\epsilon>0$, if $t_0\leq t_{C(\lambda, \epsilon)}$, we have
\begin{equation*}
P_{t_0}^{\wedge}(\mathcal{A}_C)\leq \epsilon.
\end{equation*}
In order to obtain such an upper bound, we construct a Dynkin's martingale and control its fluctuation. Let $t_0$ be a fixed time, we define a function $F\colon [0, t_0] \times \Omega_L \to \mathbb{R}$ by
\begin{equation*}
F(t, \xi)\colonequals e^{\kappa_L (t-t_0)} \Phi(\xi).
\end{equation*}
We recall that $\sigma_t^{\wedge}$, defined in Subsection \ref{graphical construction},  is the dynamics at time $t$ starting with the maximal path $\wedge$.
Further, we define a Dynkin's martingale by
\begin{equation}
M_t\colonequals F(t, \sigma^{\wedge}_t) - F(0, \sigma^{\wedge}_0) -\int_0^t (\partial_s+\mathcal{L})F(s, \sigma_s^{\wedge})ds. \label{Dynkin's martingale}
\end{equation}
 Applying $(\mathcal{L}\Phi)(\xi)=-\kappa_L \Phi(\xi)+\Psi(\xi)$ in (\ref{operator acting on weighted path}), we achieve
\begin{align}
(\partial_t+\mathcal{L})F(t, \sigma_t^{\wedge})=e^{\kappa_L (t-t_0)} \Psi(\sigma^{\wedge}_t). \label{the infinitesimal  of the weight area}
\end{align}
For simplicity of notation, set
\begin{equation}
 B(t)\colonequals\int_0^t e^{\kappa_L (s-t_0)} \Psi(\sigma^{\wedge}_s) ds. \label{notation for B(t)}
\end{equation}

Now we give an upper bound on $\mathbb{E}[ M_t^2]$ by controlling the martingale bracket $\langle M. \rangle$, which is such that the process $\big(M^2_t- \langle M. \rangle_t\big)_{t\geq 0}$ is a martingale with respect to its natural filtration. Since  there is at most one transition at each coordinate and each transition can change the value of $M_t$ in absolute value by at most $2e^{\kappa_L(t-t_0)}$, we have
\begin{equation*}
\partial_t \langle M. \rangle_t \leq \sum_{x=1}^{L-1} 4e^{2\kappa_L(t-t_0)}\leq  4Le^{2\kappa_L(t-t_0)}.
\end{equation*}
As $M_0=0$ and $\kappa_L=1-\cos\big(\frac{\pi}{L}\big)\geq \frac{\pi^2}{4L^2}$ for all $L\geq 4$, we obtain
\begin{equation}
\mathbb{E}[M_{t_0}^2]=\mathbb{E} [\langle M. \rangle_{t_0}] \leq \int_0^{t_0} 4Le^{2\kappa_L(t-t_0)} dt\leq \frac{8L^3}{\pi^2}. \label{upper bound for the second moment of the martingale}
\end{equation}
Furthermore, we give an upper bound  for the mean of  $ B(t_0)$, defined in (\ref{notation for B(t)}). 
Recalling the definitions of $\Psi$ and $\barpsi$   in (\ref{the effect of the wall}) and (\ref{absolute effect of the wall}) respectively, we have
\begin{align*}
\mathbb{E}[|B(t_0)|]&\leq \mathbb{E}\bigg[\int_0^{t_0} e^{\kappa_L(t-t_0)} \barpsi(\sigma^{\wedge}_t) dt\bigg]\\
&\leq \mathbb{E}\bigg[\int_0^{t_0} e^{\kappa_L(t-t_0)} \barpsi(\sigma^{\mu}_t) dt\bigg]\\
&\leq C_3(\lambda) \kappa_L^{-1} \sum_{x=1}^{L-1}\barsin(x) \frac{L^{3/2}}{x^{3/2}(L-x)^{3/2}}  \\
& \leq C_4(\lambda)L^{3/2}. \numberthis \label{upper bound of the absolute effect of the wall}
\end{align*}
The first inequality uses  $\vert \Psi(\xi)\vert\leq \barpsi(\xi)$ for all $\xi \in \Omega_L$. 
The second inequality is due to two facts: (1) $\barpsi(\xi)\leq \barpsi(\xi')$ for $\xi\leq \xi'$; and (2) $\sigma_t^{\wedge}\geq \sigma_t^{\mu}$. 
In the third inequality, we use Fubini's Theorem to interchange the orders of integration and expectation, and use
Lemma \ref{lemma for the maximal dynamics hitting 0 or 1} to give an upper bound for $\mathbb{E}[\barpsi(\sigma_t^{\mu})]$. In the last inequality, we use the following inequality: $$\barsin(x)=\sin \Big(\frac{\pi x}{L}\Big)\leq \frac{\min \big(x, L-x \big)\pi}{L}.$$

Here and now, we try to find the suitable small $t_0$ such that $\Phi(\sigma_{t_0}^{\wedge})$ is much larger than $L^{3/2}$ with high probability.
 We note that $\Phi(\sigma_0^{\wedge})\geq c_0 L^2$ and $\kappa_L\leq \frac{\pi^2}{2L^2}$ for all $L\geq 2$, where $c_0>0$ is a universal constant. Let $C\geq 1$, and define
$$t_0\colonequals\frac{1}{\pi^2}L^2\log L-CL^2.$$
If $t_0< 0$, nothing needs to be done (for $L\leq 4$, $t_0\leq 0$). In the remaining of this subsection, we assume $t_0\geq 0$.
Then for all $L\geq 2$, $ t_0\kappa_L\leq \frac{1}{2}\log L -C$.
Moreover, there exists a universal constant $C_0\geq 1$ such that if $C\geq C_0$, we have  $$c_0 e^{C}\geq 3C.$$
By Lemma \ref{lower bound of the area of the maximal path},  for all $C \geq \max \big(C_0, c(\lambda)\big)$, we have
\begin{equation*}
\mathbb{E}[\Phi(\sigma_{t_0}^{\wedge})]\geq 3CL^{3/2}-c(\lambda)L^{3/2}\geq 2CL^{3/2}.
\end{equation*}
Then,  if $\Phi(\sigma_{t_0}^{\wedge})\leq CL^{3/2}$ (\textit{i.e.} $\sigma_{t_0}^{\wedge}\in \mathcal{A}_C$, defined in (\ref{definition of the event AC})), it implies 
\begin{equation*}
\big \vert \Phi(\sigma_{t_0}^{\wedge})-\mathbb{E}\big[\Phi(\sigma_{t_0}^{\wedge})\big] \big\vert \geq CL^{3/2}
\end{equation*} 
and
\begin{equation}
P_{t_0}^{\wedge}(\mathcal{A}_C)\leq \mathbb{P}\big[|\Phi(\sigma_{t_0}^{\wedge})-\mathbb{E}[\Phi(\sigma_{t_0}^{\wedge})]|\geq CL^{3/2}\big]. \label{an upper bound for the maximal path dynamics  with small area}
\end{equation}
In addition, recalling $\Phi(\sigma_{t_0}^{\wedge})=F(t_0, \sigma_{t_0}^{\wedge})=M_{t_0}+F(0, \sigma_0^{\wedge})+B(t_0)$ in (\ref{Dynkin's martingale}) and using Markov's inequality,  we obtain
\begin{equation}
\begin{aligned}
& \mathbb{P}\Big[\big\vert\Phi(\sigma_{t_0}^{\wedge})-\mathbb{E}[\Phi(\sigma_{t_0}^{\wedge})]\big\vert\geq CL^{3/2}\Big]\\
=&\mathbb{P}\Big[\big\vert M_{t_0}+B(t_0)-\mathbb{E}[B(t_0)]\big\vert\geq CL^{3/2}\Big] \\
 \leq & \mathbb{P}\Big[ \big \vert M_{t_0}\big \vert\geq \frac{1}{3}CL^{3/2}\Big]+\mathbb{P}\Big[\big\vert B(t_0)\big\vert \geq \frac{1}{3}CL^{3/2}\Big]\\ 
 \leq &  \frac{9 \mathbb{E}[M_{t_0}^2]}{C^2L^3}+\frac{3\mathbb{E}[\vert B(t_0)\vert]}{CL^{3/2}}, \label{the upper bound for the well chosen event for the maximal path in the size of the windows for the lower bound of the cutoff}
 \end{aligned}
 \end{equation}
where the second last inequality holds for $C > 3C_4(\lambda)$ by  $\mathbb{E}[\vert B(t_0)\vert]\leq C_4(\lambda)L^{3/2}$ in (\ref{upper bound of the absolute effect of the wall}).
The last term in (\ref{the upper bound for the well chosen event for the maximal path in the size of the windows for the lower bound of the cutoff}) is smaller than or equal to $\epsilon/2$ for $C\geq \max \big(\frac{18}{\pi \sqrt{\epsilon}},\frac{12 C_4(\lambda)}{\epsilon}\big)$, on account of $\mathbb{E}[M_{t_0}^2]\leq \frac{8L^3}{\pi^2}$ in (\ref{upper bound for the second moment of the martingale}) and (\ref{upper bound of the absolute effect of the wall}). 
Combining (\ref{lower bound for the well chosen event of the equilibrium states for the size of the windows}), (\ref{an upper bound for the maximal path dynamics  with small area}) and (\ref{the upper bound for the well chosen event for the maximal path in the size of the windows for the lower bound of the cutoff}), we know that
\begin{equation}
\Vert P^{\wedge}_{t_0}-\mu \Vert_{\TV}\geq \mu(A_C)-P_{t_0}^{\wedge}(A_C)\geq 1-\epsilon, \label{lower bound for the two events}
\end{equation}
which holds for $C\geq \max\{\frac{2c}{\epsilon}, C_0, 3C_4(\lambda), \frac{18}{\pi \sqrt{\epsilon}},\frac{12C_4(\lambda)}{\epsilon} \} \equalscolon C(\lambda, \epsilon)$. 
Therefore, for $C\geq C(\lambda, \epsilon)$,  we have 
\begin{equation*}
\T(\epsilon)\geq \frac{1}{\pi^2}L^2\log L-CL^2.
\end{equation*}
\end{proof}
\section{Upper bound on the mixing time for $\lambda \in (0, 1]$}\label{upper bound for the mixing time in the first diffusive regime}
This section is devoted to providing an upper bound on the mixing time of the dynamics for the regime $\lambda \in (0, 1]$.
For any $\xi\in \Omega_L$, by the  triangle inequality, we have
\begin{equation}
\Vert P_t^{\xi}-P^{\mu}_t \Vert_{\TV}\leq \sum_{\xi' \in \Omega_L}\mu(\xi') \Vert P_t^{\xi}-P_t^{\xi'} \Vert_{\TV}\leq \max_{\xi'\in \Omega_L} \Vert P^{\xi}_t-P^{\xi'}_t\Vert_{\TV}. \label{preparation inequality1 for using the coupling time}
\end{equation}
To give an upper bound for the term in the rightmost hand side above,
we use the following characterization of total variation distance.
Let $\alpha$ and $\beta$ be two probability measures on $\Omega_L$. We say that $\vartheta$ is a coupling of $\alpha$ and $\beta$, if $\vartheta$ is a probability measure on $\Omega_L \times \Omega_L$ such that $\vartheta(\xi \times \Omega)=\alpha(\xi)$ and $\vartheta(\Omega\times \xi')=\beta(\xi')$ for any elements $\xi,\xi' \in \Omega_L$.  The following proposition says that the total variation distance measures how well we can couple two random variables with distribution laws $\alpha$ and $\beta$ respectively.
\begin{proposition}[Proposition 4.7 \cite{LPWMCMT}] \label{characterization of total variation distance}
Let $\alpha$ and $\beta$ be two probability distributions on $\Omega_L$. Then
\begin{equation*}
\Vert\alpha -\beta\Vert_{\TV}=\inf \big\{ \vartheta(\xi \neq \xi'):   \mbox{ $\vartheta$ is a coupling of  $\alpha$ and  $\beta$}  \big\}.
\end{equation*}
\end{proposition}

The graphical construction in Subsection \ref{graphical construction} provides a coupling between $P_t^{\xi}$ and $P_t^{\xi'}$, which preserves the monotonicity asserted in Proposition \ref{preserving the monotonicity}. Therefore, $\sigma_t^{\xi}$ lies between $\sigma_t^{\vee}$ and $\sigma_t^{\wedge}$ for any $\xi\in \Omega_L$.
Applying Proposition \ref{characterization of total variation distance}, we obtain
\begin{equation}
\Vert P^{\xi}_t -P_t^{\xi'}\Vert_{\TV}\leq \mathbb{P}[\sigma_t^{\xi}\neq \sigma_t^{\xi'}]\leq \mathbb{P}[\sigma_t^{\wedge}\neq \sigma_t^{\vee}], \label{preparation inequality2 for using the coupling time}
\end{equation}
where the last inequality is due to the fact that after the dynamics starting from the two extremal paths have coalesced, we must have $\sigma_t^{\wedge}=\sigma_t^{\xi}=\sigma_t^{\vee}$  for any $\xi\in \Omega_L$.  This argument was used in \cite[Theorem 3.1]{caputo2008approach} to obtain an upper bound on the mixing time. 
Comparing with the coupling in \cite[Subsection 2.2.1]{caputo2008approach},
the graphical construction in Subsection \ref{graphical construction} provides more independent flippable corners and maximizes the fluctuation of the area enclosed by $\sigma_t^{\wedge}$ and $\sigma_t^{\vee}$.
Adapting the approach in \cite[Section 7]{labbe2018mixing}, we use a supermartingale approach to have a good control of the fluctuation of the area enclosed by $\sigma_t^{\wedge}$ and $\sigma_t^{\vee}$ to achieve a sharp upper bound on the mixing time.
Let the coalescing time $\tau$ be
\begin{equation*}
 \tau\colonequals\inf \{t\geq 0: \sigma_t^{\wedge}=\sigma_t^{\vee} \},  
\end{equation*} 
which is the first moment when the dynamics starting from the two extremal paths coalesce.
By (\ref{preparation inequality1 for using the coupling time}) and (\ref{preparation inequality2 for using the coupling time}), we obtain
\begin{equation}
d^{L, \lambda}(t)\leq  \mathbb{P}[\sigma_t^{\wedge}\neq \sigma_t^{\vee}]=\mathbb{P}[\tau>t]. \label{upper bound for mixing time in terms of coupling time}
\end{equation}
In this section, our goal is to show that for any given $\delta>0$ and all $L$ sufficiently large, with high probability,  we have
\begin{equation*}
 \tau \leq  \frac{1+\delta}{\pi^2} L^2\log L.
\end{equation*}
We adapt the approach in \cite[Section 7]{labbe2018mixing} to achieve this goal.
In practice, it is more feasible to couple  two dynamics when, at least, one of them is at equilibrium. 
Let
\begin{equation}
\begin{aligned}
\tau_1:&=\inf \{ t\geq 0, \sigma^{\wedge}_t=\sigma^{\mu}_t  \}, \\ \label{absorption time for A}
\tau_2:&=\inf \{ t\geq 0, \sigma^{\vee}_t=\sigma^{\mu}_t  \},
\end{aligned}
\end{equation}
where we recall that the dynamics $(\sigma_t^{\mu})_{t\geq 0}$ starts with the equilibrium distribution $\mu$, 
defined in Subsection \ref{graphical construction}.
By the definition of $\tau$, we know that
\begin{equation*}
\tau= \max \big(\tau_1, \tau_2 \big).
\end{equation*}
For this goal, it is sufficient to prove the following proposition.
\begin{proposition}\label{proposition of merging time for the extremal paths}
For $i\in \{1,2\}$, any given $\lambda\in (0, 1]$ and  $\delta>0$, we have
\begin{eqnarray}
\lim_{L \rightarrow \infty}\mathbb{P}\Big[\tau_i \leq (1+\delta)\frac{1}{\pi^2}L^2\log L\Big]=1. \label{merging time of the maximal path}
\end{eqnarray}
\end{proposition}
Theorem \ref{the cutoff phenomenon in the diffusive regime} is proved as a combination of Proposition \ref{lower bound of the mixing time} and Proposition \ref{proposition of merging time for the extremal paths}. Therefore, there is a cutoff in the Markov chains for $\lambda\in (0, 1]$.
Since the proofs about $\tau_1$ and $\tau_2$ in Proposition \ref{proposition of merging time for the extremal paths}  are similar,
we  only give the proof of (\ref{merging time of the maximal path}) for $\tau_1$.
For any given $\delta>0$,  set $$t_{\delta}\colonequals(1+\delta)\frac{1}{\pi^2}L^2\log L.$$
We outline the idea for the proof.  We define a weighted area function $A_t$, in (\ref{definition of A_t}) below, which is almost the area enclosed by the paths $\sigma_t^{\wedge}$ and $\sigma_t^{\mu}$ at time $t$. By contraction property of the area function $A_t$, $\sigma_t^{\wedge}$ gets fairly close to the path $\sigma_t^{\mu}$ at time $t_{\delta/2}$. After time $t_{\delta/2}$, 
we estimates the fluctuation of function $A_t$ by supermartingale approach. Then we compare the fluctuation of the function $A_t$ with intermediate time increments to obtain (\ref{merging time of the maximal path}).

\subsection{A weighted area function.}  In this subsection, we define an area function $A_t$ and base on $A_t$ to study  $\tau_1$.
 First, inspired by \cite[Equation (1)]{wilson2004mixing}, we define a function $\barphi_{\beta}\colon \Omega_L  \to [0, \infty)$ given by
 \begin{equation*}
    \barphi_{\beta}(\xi)\colonequals\sum_{x=1}^{L-1} \xi_x \barcos_{\beta}(x),
 \end{equation*}
where $\barcos_{\beta}(x)\colonequals\cos\Big(\frac{\beta (x-L/2)}{L}\Big)$, and $\beta$ is a constant in $(2 \pi/3, \pi)$. The constant $\beta$ is only dependent on $\delta$ and sufficiently close to $\pi$, which will be chosen
 in the proof of Lemma \ref{the first stopping time  behaves well} below.
We can see that $\barphi_{\beta}(\xi)$ is approximately the area enclosed by the $x$-axis and the path $\xi \in \Omega_L$. 
  Throughout this paper, we omit the index $\beta$ in $\barphi_{\beta}$ and $\barcos_{\beta}$ as much as possible. Observe that if $\xi$ and $\xi'$ are two elements of $\Omega_L$ satisfying $\xi\leq \xi'$, then
 \begin{equation}
 \barphi(\xi)\leq \barphi(\xi'). \label{monotonicity of the function barphi}
 \end{equation}
  The minimal increment of the function $\barphi$ is
\begin{align}
\delta_{\min}&\colonequals\min_{\xi\leq \xi', \xi\neq \xi'} \Big(\barphi(\xi')-\barphi(\xi)\Big)
=2\cos\Big(\frac{\beta(L/2-1)}{L}\Big), \label{minimal increment of A_t}
\end{align}
and 
\begin{equation}
2\cos\Big(\frac{\beta(L/2-1)}{L}\Big)
 \geq \frac{1}{2} (\pi-\beta) 
\end{equation}
for $L\geq 6$ and $\beta\in (2\pi/3, \pi)$, where we use the inequality $\cos(\pi/2-x)=\sin x \geq x/2$ for $x \in [0, \pi/3]$.
 Let the weighted area function $A\colon [0, \infty) \to [0, \infty)$ be
\begin{equation}
A_t\colonequals\frac{\barphi(\sigma_t^{\wedge})-\barphi(\sigma_t^{\mu})}{\delta_{\min}}. \label{definition of A_t}
\end{equation}
We observe that $\tau_1$, defined in (\ref{absorption time for A}), is the first time  at which  $A_t$ reaches zero. Moreover, $A_t$ equals to zero if and only if $\sigma_t^{\wedge}$ equals to $\sigma_t^{\mu}$.  If $\tau_1\leq t_{\delta/2}$, we are done.  In the rest of this section, we assume $\tau_1>t_{\delta/2}$.

Take $\eta>0$ and sufficiently small, and  $K\colonequals\lceil 1/(2\eta) \rceil$.
 We define a sequence of successive stopping times $(\Tau_i)_{i=2}^{K}$ by
\begin{equation*}
\Tau_2 \colonequals \inf \Big\{ t\geq t_{\delta/2} :  A_t\leq L^{\frac{3}{2}-2\eta}  \Big\},
\end{equation*}
and for  $3\leq i \leq K$,
\begin{equation*}
 \Tau_i \colonequals \inf \Big\{ t\geq \Tau_{i-1}: A_t\leq L^{\frac{3}{2}-i\eta}   \Big\} .
\end{equation*}
 For consistency of notations, we set $\Tau_{\infty}\colonequals\max \big(\tau_1, t_{\delta/2}\big)$.  
The remaining of this section is devoted to proving the following proposition.
\begin{proposition}\label{lemma for controlling the increment for consecutive stopping times}  Given $\delta>0$, if $\eta$ is chosen to be a sufficiently small positive constant with $K \colonequals\lceil 1/(2\eta) \rceil>1/(2\eta)$, we have
\begin{equation*}
\lim_{L\rightarrow \infty}\mathbb{P} \bigg[ \{\Tau_2=t_{\delta/2}\} \cap  \bigg(  \bigcap_{i=3}^{K} \{\Delta \Tau_i \leq 2^{-i}L^2 \}\bigg)  \cap \{\Tau_{\infty}-\Tau_K \leq L^2\} \bigg]=1, 
\end{equation*}
where $\Delta \Tau_i\colonequals\Tau_i -\Tau_{i-1}$ for $3 \leq i \leq K$.
\end{proposition}
\noindent
If Proposition \ref{lemma for controlling the increment for consecutive stopping times} holds,  for $L$ sufficiently large, we have
\begin{equation*}
\tau_1=\Tau_{\infty}\leq t_{\delta/2}+\sum_{i=3}^{K}2^{-i}L^2+L^2\leq (1+\delta)\frac{1}{\pi^2}L^2 \log L.
\end{equation*}
Hence, Proposition \ref{proposition of merging time for the extremal paths}  is proved.
We sketch the idea for Proposition \ref{lemma for controlling the increment for consecutive stopping times} as follows. Due to contraction property of the function $A_t$, we obtain $\Tau_2=t_{\delta/2}$. In addition, we compare the increments $\Delta \Tau_i$, $\Tau_{\infty}-\Tau_K$ with the fluctuation of function $A_t$.

\subsection{The proof of $\Tau_2=t_{\delta/2}$.} The main task of this subsection is to prove that  the function $A_t$ has a contraction property, due to which we obtain $\Tau_2=t_{\delta/2}$ with high probability.
Above all, we want to understand how  the generator $\mathcal{L}$ acts on the function $\barphi$. We have
\begin{equation*}
(\mathcal{L}\barphi)(\xi)=\sum_{x=1}^{L-1}\barcos(x) \mathcal{L}\xi_x.
\end{equation*}
We recall Lemma \ref{the action of operator on paths}: for any $\xi\in \Omega_L$,
\begin{equation*}
\mathcal{L}\xi_x=(\Delta \xi)_x +\mathbbm{1}_{\{\xi_{x-1}=\xi_{x+1}=0 \}}+\Big(\frac{1-\lambda}{1+\lambda}\Big)\mathbbm{1}_{\{\xi_{x-1}=\xi_{x+1}=1\}}.
\end{equation*}
For $\xi, \xi'\in \Omega_L$, we have
\begin{equation}
\sum_{x=1}^{L-1}\barcos(x)\Big((\Delta \xi')_x-(\Delta \xi)_x\Big)= -\Big(1-\cos\Big(\frac{\beta}{L}\Big)\Big)\sum_{x=1}^{L-1}\barcos(x)(\xi'_x-\xi_x). \label{Monotonicity of the new weighted area}
\end{equation}
Considering
\begin{equation*}
\mathcal{L}\xi_x-(\Delta \xi)_x=\mathbbm{1}_{\{\xi_{x-1}=\xi_{x+1}=0 \}}+\Big(\frac{1-\lambda}{1+\lambda}\Big)\mathbbm{1}_{\{\xi_{x-1}=\xi_{x+1}=1\}},
\end{equation*}
we see that both terms in the right-hand side are nonnegative and monotonically decreasing in $\xi$ for $\lambda \in (0, 1]$.
Hence, if $\xi\leq \xi'$, we know that
\begin{equation}
\mathcal{L}\xi_x-(\Delta \xi)_x \geq \mathcal{L}\xi'_x-(\Delta \xi')_x. \label{Monotonicity of the effect of the wall}
\end{equation} 
For simplicity of notation, we set 
\begin{equation*}
\bargamma=\bargamma_{L, \beta}\colonequals1-\cos(\beta/L).
\end{equation*}
 On the grounds of Lemma \ref{the action of operator on paths}, (\ref{Monotonicity of the new weighted area}) and (\ref{Monotonicity of the effect of the wall}), if $\xi\leq \xi'$, we obtain
\begin{align*}
(\mathcal{L}\barphi)(\xi')-(\mathcal{L}\barphi)(\xi)&=\sum_{x=1}^{L-1}\barcos(x) \Big((\Delta \xi')_x -(\Delta \xi)_x+(\mathcal{L}\xi'_x-(\Delta \xi')_x) -(\mathcal{L}\xi_x-(\Delta \xi)_x )\Big)\\
&\leq \sum_{x=1}^{L-1}\barcos(x)\Big((\Delta \xi)_x' -(\Delta \xi)_x\Big)\\
&\leq -\bargamma \Big(\barphi(\xi')-\barphi(\xi)\Big). \numberthis \label{preparation for the lemma about contraction for the first time}
\end{align*}
Now we are ready to prove that $\Tau_2=t_{\delta/2}$ with high probability.
\begin{lemma}\label{the first stopping time  behaves well} For all  $\epsilon>0$, all sufficiently small $\delta>0$ and $0<\eta<\delta/10$,  if  $L$ is sufficiently large, we have  
\begin{equation*}
\mathbb{P} \big[\Tau_2>t_{\delta/2}\big]\leq \epsilon.
\end{equation*}
\end{lemma}
\begin{proof}
By $\sigma_t^{\wedge}\geq \sigma_t^{\mu}$ and (\ref{preparation for the lemma about contraction for the first time}), we obtain
\begin{equation}
\begin{aligned}
\frac{d}{dt}\mathbb{E}\Big[\barphi(\sigma_t^{\wedge})-\barphi(\sigma_t^{\mu})\Big]&=\mathbb{E}\Big[(\mathcal{L}\barphi)(\sigma_t^{\wedge})-(\mathcal{L}\barphi)(\sigma_t^{\mu})\Big]\\
&\leq -\bargamma \mathbb{E}\Big[\barphi(\sigma_t^{\wedge})-\barphi(\sigma_t^{\mu})\Big]. \label{the derivative of the new weighted area} 
\end{aligned}
\end{equation}
Using (\ref{the derivative of the new weighted area}),  $\barphi(\sigma_0^{\wedge})\leq \frac{1}{2}L^2$ and $\barphi(\xi)\geq 0$  for all $\xi \in \Omega_L$, we obtain
\begin{equation}
\mathbb{E}\Big[\barphi(\sigma_t^{\wedge})-\barphi(\sigma_t^{\mu})\Big]\leq e^{-\bargamma t}\Big(\barphi(\sigma_0^{\wedge})-\barphi(\sigma_0^{\mu})\Big)\leq \frac{1}{2} L^2 e^{-\bargamma t}. \label{upper bound for the expectation of the difference new weighted are between maximal path and equilibrium path}
\end{equation}
Thus, applying Markov's inequality, we achieve
\begin{equation}
\begin{aligned}
\mathbb{P}[\Tau_2>t_{\delta/2}]&=\mathbb{P}[A_{t_{\delta/2}}>L^{\frac{3}{2}-2\eta}]\\
&\leq \frac{1}{2\delta_{\min}}L^{2\eta+\frac{1}{2}}e^{-\bargamma t_{\delta/2}},
\label{the upper bound of the probability of the event that the weighted  area doesn't shrink}
\end{aligned}
\end{equation}
where the last inequality uses (\ref{upper bound for the expectation of the difference new weighted are between maximal path and equilibrium path}) and  the definition of $A_t$ in (\ref{definition of A_t}).
Recalling $\bargamma=1-\cos(\beta/L)$ and using the inequality $1-\cos x \geq \frac{1}{2}x^2-\frac{1}{24}x^4$ for all $x \geq 0$, we have
$$\bargamma t_{\delta/2}\geq \frac{\beta^2}{2\pi^2}(1+\frac{\delta}{2})\log L-\frac{\beta^4}{24 L^2}(1+\frac{\delta}{2}) \log L.$$ 
For $\delta>0$ sufficiently small and $0<\eta<\delta/10$, we choose
$$\beta=\pi \sqrt{\frac{1+\frac{9}{20}\delta}{1+\frac{\delta}{2}}}\in (2\pi/3, \pi)$$
which satisfies
\begin{equation*}
\frac{1}{2}(1+\frac{\delta}{2})\frac{\beta^2}{\pi^2}=\frac{1}{2}+\frac{9}{40}\delta > \frac{1}{2}+2\eta.
\end{equation*}
 With this choice of $\beta$, the rightmost term of (\ref{the upper bound of the probability of the event that the weighted  area doesn't shrink}) vanishes as $L$ tends to infinity.
\end{proof}

\subsection{ The estimation of  $\langle A. \rangle_{\Tau_i}-\langle A. \rangle_{\Tau_{i-1}}$.} \label{diffusion estimate after contraction} 
Due to Dynkin's martingale formula, we know that
\begin{equation*}
A_t-A_0-\int_{0}^t \mathcal{L}A_s ds
\end{equation*}
is a martingale. Moreover, we let $\langle A. \rangle_t$ represent the predictable bracket associated with this martingale. 
The objective of this subsection is to show that $\langle A. \rangle_{\Tau_i}-\langle A. \rangle_{\Tau_{i-1}}$ is small for all $i\in \lint 3, K \rint$.
For any $i \in \lint 3, K \rint$, let
\begin{equation}
\begin{aligned}
&\Delta_i \langle A \rangle\colonequals\langle A. \rangle_{\Tau_i}- \langle A. \rangle_{\Tau_{i-1}},
\end{aligned}
\end{equation}
and let
\begin{equation}
\Delta_{\infty} \langle A \rangle\colonequals \langle A. \rangle_{\Tau_{\infty}}- \langle A. \rangle_{\Tau_K}.
\end{equation}
We have $\mathcal{L}A_s\leq 0$, according to  (\ref{preparation for the lemma about contraction for the first time}), $\sigma_t^{\wedge}\geq \sigma_t^{\mu}$, and the monotonicity of the function $\barphi$ stated in (\ref{monotonicity of the function barphi}).
 Then, $A_t$ is a supermartingale for $\lambda \in (0, 1]$.  
Its  jump amplitudes in absolute value are bounded below by $1$ before the absorption time  $\tau_1$, which is defined in (\ref{absorption time for A}). Additionally, by  the graphical construction in Subsection \ref{graphical construction}, the jump rates of $A_t$ are least $1$ before the absorption time $\tau_1$, 
because
there are at least two flippable corners in $\sigma_t^{\wedge}$ or $\sigma_t^{\mu}$ that can change the value of $A_t$ for $t< \tau_1$. We refer to Figure \ref{fig:jumprates} in Subsection \ref{Preliminaries} for illustration:  the square formed by $\xi$ (\textit{i.e.} the blue path) and $\xi^x$ (\textit{i.e.} the blue path with a black dashed corner) are the flippable corners that can change the value of $A_t$, and their rates in total are $1$.
Now, we are in the setting to apply  
 \cite[Proposition 29]{labbe2018mixing}  which, under some condition, allows to control hitting times of supermartingales in terms of the martingale bracket.

\begin{proposition}\cite[Proposition 29]{labbe2018mixing} \label{supermartingale proposition}
Let $(\mathbf{M}_t)_{t\geq 0}$ be a pure-jump supermartingale with bounded  jump rates and jump amplitudes, and $\mathbf{M}_0\leq a$ almost surely. Let $\langle \mathbf{M}. \rangle$, with an abuse of notation, denote the predictable bracket associated with the martingale $\overline{\mathbf{M}}_t=\mathbf{M}_t-I_t$ where $I$ is the compensator of $\mathbf{M}$. 
\begin{itemize}
\item[(i)] Set $\tau_0 \colonequals\inf \{t \geq 0 :  \mathbf{M}_t=0 \}$. Assume that $\mathbf{M}_t$
is nonnegative and that,
until the absorption time $\tau_0$, its jump amplitudes and jump rates are bounded
below by $1$. For any $u> 0$,  if $a\geq 1$, we have
\begin{equation*}
\mathbb{P}\big[\tau_0 \geq a^2 u\big]\leq 4u^{-1/2}.
\end{equation*}
\item[(ii)]  Given $b \in \mathbb{R}$ and $b\leq a$, we set $ \tau_b \colonequals\inf\{t \geq 0 : \mathbf{M}_t\leq b\}$. If the
amplitudes of the jumps of $(\mathbf{M}_t)_{t\geq 0}$ are bounded above by $a-b$, for
any $u\geq 0$, we have
\begin{equation*}
\mathbb{P}\big[\langle \mathbf{M}.\rangle_{\tau_b}\geq (a-b)^2u \big]\leq 8u^{-1/2}.
\end{equation*}
\end{itemize}
\end{proposition}

 \noindent Now we apply Proposition \ref{supermartingale proposition} to prove that the event
\begin{equation*}
\mathcal{A}_L\colonequals\Big\{ \forall i\in \lint 3, K \rint,  \Delta_i \langle A \rangle \leq L^{3-2(i-1)\eta+\frac{1}{2}\eta}\Big\} \bigcap \Big\{ \Delta_{\infty} \langle A \rangle\leq L^2 \Big\}
\end{equation*}
 has almost the full mass, which is the following lemma.
\begin{lemma} \label{the event AL has almost the full mass} 
We have
\begin{equation}
\lim_{L\rightarrow \infty} \mathbb{P}[  \mathcal{A}_L]=1.
\end{equation}
\end{lemma}
\begin{proof}
We just need to show that the probability of its complement $\mathcal{A}_L^{\complement}$ is almost zero.
We apply Proposition \ref{supermartingale proposition}-(ii) to $(A_{t+\Tau_{i-1}})_{t\geq 0}$ with $a_i=L^{\frac{3}{2}-(i-1)\eta}$ and $b_i=L^{\frac{3}{2}-i\eta}$.  For every $ i \in \lint 3, K\rint$,  we obtain 
\begin{equation}
\mathbb{P}\Big[\Delta_i \langle A \rangle \geq (L^{\frac{3}{2}-(i-1)\eta}-L^{\frac{3}{2}-i\eta})^2u_i \Big ]
\leq 8 u_i^{-\frac{1}{2}}, \label{application of the supermartingale inequality}
\end{equation}
where we choose $u_i=L^{\frac{1}{2}\eta}(1-L^{-\eta})^{-2}$, satisfying
 $$\Big(L^{\frac{3}{2}-(i-1)\eta}-L^{\frac{3}{2}-i\eta}\Big)^2u_i=L^{3-2(i-1)\eta+\frac{1}{2}\eta}.$$ 
 We see that $u_i$ tends to infinity as $L$ tends to infinity. Accordingly, the rightmost term in (\ref{application of the supermartingale inequality}) vanishes as $L$ tends to infinity.

 We apply Proposition \ref{supermartingale proposition}-(i) to $(A_{t+\Tau_K})_{t\geq 0}$ with $a_{\infty}=L^{\frac{3}{2}-K\eta}$ and $b_{\infty}=0$. We choose $u_{\infty}$ such that $(a_{\infty}-b_{\infty})^2u_{\infty}=L^2$, \textit{i.e.} $$u_{\infty}=L^{-1+2K\eta},$$ which tends to infinity due to $K=\lceil 1/(2\eta) \rceil>1/(2\eta)$. Thus
$\mathbb{P}[\Delta_{\infty} \langle A \rangle \geq L^2]$ tends to zero as $L$ tends to infinity.
Since $K$ is a constant, we have
\begin{equation*}
\lim_{L\rightarrow \infty} \mathbb{P} \Big[ \mathcal{A}_L^{\complement} \Big]=0.
\end{equation*}
\end{proof}
\subsection{The comparison of $\Tau_i-\Tau_{i-1}$ to  $\Delta_i \langle A \rangle$.} As explained in Subsection \ref{diffusion estimate after contraction}, we have $\partial_t \langle A. \rangle\geq 1$ for all $t< \Tau_{\infty}$.
 Therefore,  we obtain
\begin{equation*}
\Delta_{\infty} \langle A \rangle=\int_{\Tau_K}^{\Tau_{\infty}} \partial_t \langle A. \rangle dt\geq \int_{\Tau_K}^{\Tau_{\infty}} 1  dt=\Tau_{\infty}-\Tau_K.
\end{equation*}
Hence, when the event $\mathcal{A}_L$ holds, we obtain 
\begin{equation*}
\Tau_{\infty}-\Tau_K\leq \Delta_{\infty}\langle  A \rangle \leq L^2.
\end{equation*}

Now we control the intermediate increment $\Tau_i-\Tau_{i-1}$ for $ 3 \leq i\leq K$. To do that,
we compare $\Tau_i-\Tau_{i-1}$ with $\langle A. \rangle_{\Tau_i}-\langle A. \rangle_{\Tau_{i-1}}=\Delta_i \langle A. \rangle $.
 First,  we give a lower bound on $\partial_t \langle A. \rangle$, which is related with: (a) 
the maximal contribution among all the coordinates $x\in \lint 0, L\rint$ in the definition of $A_t$; 
 and (b) the amount of flippable corners in $\sigma_t^{\mu}$ or $\sigma_t^{\wedge}$ that can change the value of $A_t$.
 Considering the definition of $A_t$ in  (\ref{definition of A_t}),
set 
\begin{equation}
H(t)\colonequals\max_{x \in \lint 0, L \rint} \sigma_t^{\wedge}(x). \label{definitionofH}
\end{equation}
For a lower bound on the quantity mentioned in (b), we need the maximal length of the monotone segments of $\sigma_t^{\mu}$.  
For $\xi \in \Omega_L$, we define
\begin{align*}
Q_1(\xi)&\colonequals\max \{n\geq 1, \exists i \in \lint 0, L-n\rint, \forall x \in \lint i+1, i+n \rint, \xi_x-\xi_{x-1}=1\},\\
Q_2(\xi)& \colonequals\max\{n\geq 1, \exists i \in \lint 0, L-n \rint, \forall x \in \lint  i+1, i+n\rint, \xi_x-\xi_{x-1}=-1\},
\end{align*}
and
\begin{equation}
 Q(\xi)\colonequals\max \big(Q_1(\xi), Q_2(\xi)\big).
\label{maximal length of monotone segment}
\end{equation}
Using these two quantities $H(t)$ and $Q(\sigma_t^{\mu})$, we obtain a lower bound for $\partial_t \langle A. \rangle$, which is the following lemma.
\begin{lemma}We have \label{lemma for the lower bound of the derivative of the second weighted area}
\begin{equation}
\partial_t\langle A. \rangle\geq \max \bigg( 1, \frac{\lambda \delta_{\min} A_t}{3(1+\lambda) H(t) Q(\sigma^{\mu}_t)} \bigg).
\end{equation}
\end{lemma}
\begin{proof}
We observe that $A_t$ displays a jump whenever either $\sigma^{\mu}_t$ or $\sigma_t^{\wedge}$  flips a corner.
Note that by (\ref{definition of A_t}) and (\ref{minimal increment of A_t}),
 any jump amplitude in absolute value of $A$ is at least $1$.
Since any flippable corner is flipped with rate at least $$\min \Big\{ \frac{1}{2}, \frac{1}{1+\lambda}, \frac{\lambda}{1+\lambda}\Big\}=\frac{\lambda}{1+\lambda},$$ we obtain 
\begin{equation*}
\partial_t \langle A.\rangle_t\geq \frac{\lambda}{1+\lambda} \# \big\{x\in \mathcal{B}_t: \Delta \sigma^{\mu}_t(x)\neq 0  \big\}
\end{equation*}
where 
$\mathcal{B}_t\colonequals\big\{  x\in \lint 1, L-1 \rint : \exists y \in \lint x-1, x+1 \rint,  \sigma_t^{\wedge}(y) \neq  \sigma_t^{\mu}(y) \big\}.$ For simplicity of notation, set  $\mathcal{D}_t \colonequals \big\{x\in \mathcal{B}_t: \Delta \sigma^{\mu}_t(x)\neq 0  \big\} $.
Let $\lint a, b \rint$ denote the horizontal coordinates of a maximal connected component of $\mathcal{B}_t$, for which we refer to Figure \ref{fig:lowerboundforderivativeofA} for illustration. Since $\sigma_t^{\mu}$ can  not be monotone in the entire domain $\lint a, b \rint$, we know that
\begin{equation*}
\# (\mathcal{D}_t \cap  \lint a, b \rint) \geq 1.
\end{equation*}

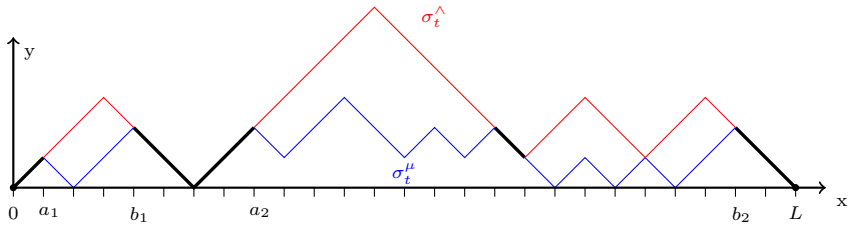
\begin{figure}[h]
\centering
  \begin{tikzpicture}[scale=.4,font=\tiny]
   \draw (25,4) -- (25,-1) -- (52,-1);
   \draw[color=black, line width=0.45mm](25,-1)--(26,0); \draw[color=black, line width=0.45mm](29,1)--(30,0)--(31,-1)--(33,1);
    \draw[color=black, line width=0.45mm](41,1)--(42,0);
   \draw[color=black, line width=0.45mm](49,1) -- (50,0)--(51, -1);
    \draw[color=blue] (26,0) -- (27,-1) -- (28,0) -- (29,1);    
      \draw[color=blue]  (33,1)-- (34,0) -- (35,1) -- (36,2) -- (37,1) -- (38,0) -- (39,1) -- (40,0) -- (41,1); 
      \draw[color=blue] (42,0) -- (43,-1) -- (44,0) -- (45,-1) -- (46,0)--(47,-1)--(48,0)--(49,1);
   \draw[color=red](26,0)--(27,1)--(28,2)--(29,1);     
   \draw[color=red](33,1)--(37,5)--(41,1);
    \draw[color=red](42,0)--(44,2)--(46,0)-- (47,1) -- (48,2) -- (49,1); 
     \node[below] at (26.2,-1.3) {$a_1$};
      \node[below] at (29.2,-1.3) {$b_1$};
     \node[below] at (33.2,-1.3) {$a_2$};
     \node[below] at (49.2,-1.3) {$b_2$};
    \foreach \x in {25,...,51} {\draw (\x,-1.3) -- (\x,-1);}
    \draw[fill] (25,-1) circle [radius=0.1];
    \draw[fill] (51,-1) circle [radius=0.1];   
    \node[below] at (25,-1.3) {$0$};
    \node[below] at (51,-1.3) {$L$};
    \node[red,above] at (39, 4){$\sigma_t^{\wedge}$};
    \node[blue,above] at (38, -1.2){$\sigma_t^{\mu}$};
    \draw[thick,->] (25,-1) -- (25,4) node[anchor=north west]{y};
    \draw[thick,->] (25,-1) -- (52,-1) node[anchor=north west]{x};
  \end{tikzpicture}
  \caption{\label{fig:lowerboundforderivativeofA} In this figure, $\sigma_t^{\wedge}$ consists of the red line segments and black thick line segments, while $\sigma_t^{\mu}$ consists of the blue line segments and black thick line segments. Moreover, $\mathcal{B}_t= \lint  a_1, b_1\rint \cup \lint a_2, b_2\rint$, $\#(\mathcal{D}_t \cap \lint a_1, b_1\rint)=3$, and $\#(\mathcal{D}_t \cap \lint a_2, b_2\rint)=13$. In $\lint a_2, b_2\rint$, the monotone segments of $\sigma_t^{\mu}$ are $\lint a_2, a_2+1\rint$, $\lint a_2+1, a_2+3\rint$, $\lint a_2+3, a_2+5\rint$, and so on as shown in the figure.}
\end{figure}

\noindent
In $\mathcal{B}_t$, we decompose the path associated with $\sigma^{\mu}_t$ into consecutive maximal  monotone segments. Then we know that in $\mathcal{B}_t$ every two consecutive components correspond to one flippable corner, which is a point in $\mathcal{D}_t$. As any maximal monotone component is at most of length  $Q(\sigma_t^{\mu})$ defined in (\ref{maximal length of monotone segment}), 
we obtain
\begin{equation} 
\# (\mathcal{D}_t \cap \lint a, b \rint) \geq \frac{1}{2} \bigg\lfloor \frac{b-a}{Q(\sigma_t^{\mu})} \bigg\rfloor \geq \frac{1}{3} \frac{b-a}{Q(\sigma_t^{\mu})}.\label{lower bound for the number of flippable corners}
\end{equation}
In addition, we observe that 
\begin{equation} 
\sum_{x=a}^{b}\frac{(\sigma_t^{\wedge}(x)-\sigma_t^{\mu}(x)) \barcos(x)}{\delta_{\min}}\leq (b-a) \frac{H(t)}{\delta_{\min}}, \label{upper bound for the weighted area by using the highest height}
\end{equation}
where $H(t)$ is defined in (\ref{definitionofH}).
Summing up all such intervals $\lint a, b \rint$ and using  (\ref{lower bound for the number of flippable corners}) and (\ref{upper bound for the weighted area by using the highest height}), we obtian
\begin{equation*}
A_t\leq \frac{3}{\delta_{\min}} H(t) Q(\sigma_t^{\mu}) \#\mathcal{D}_t.
\end{equation*} 
Therefore, we have 
\begin{equation*}
\partial_t \langle A. \rangle \geq \frac{\lambda }{1+\lambda} \# \mathcal{D}_t \geq   \frac{\lambda \delta_{\min}}{3(1+\lambda)} \frac{A_t}{H(t) Q(\sigma_t^{\mu})}.
\end{equation*}
This yields the desired result.
\end{proof}

 To give a good lower bound for $\partial_t \langle A. \rangle $,   we need to  control $Q(\sigma_t^{\mu})$ and $ H(t)$. Our next step is to give an upper bound on $Q(\sigma_t^{\mu})$, which is the following lemma. We recall the notation $$t_{\delta}=(1+\delta)\frac{1}{\pi^2}L^2\log L.$$
 \begin{lemma}\label{flippable corners for the equilibrium} We have
\begin{equation}
\lim_{L\rightarrow \infty} \mathbb{P} \Big[\exists t \in [0, t_{\delta}]: Q(\sigma_t^{\mu})> (\log L)^2\Big]=0. \label{upper bound for flippable corners}
\end{equation}
\end{lemma}
\begin{proof}
Firstly, we prove that there exists a constant $C(\lambda)>0$ such that for all $L\geq 2$
\begin{equation}
\mu(Q(\xi)>(\log L)^2)\leq 2C(\lambda)L^{5/2}2^{-(\log L)^2}. \label{the estimation in the equilibrium state for flippable corners}
\end{equation}
Since there are at most $L$ starting positions for a monotone segments either monotonically increasing or decreasing, we have $$\# \{\xi\in \Omega_L: Q(\xi)>(\log L)^2 \}\leq L2^{1+L-(\log L)^2}.$$ Moreover,  as $\lambda^{\mathcal{N(\xi)}}\leq 1$ for $\lambda \in (0, 1]$ and any $\xi\in \Omega_L$, we obtain
\begin{equation}
\mu(Q(\xi)>(\log L)^2)\leq C_5(\lambda)\frac{L2^{1+L-(\log L)^2}}{2^L L^{-3/2}}=2C_5(\lambda)L^{5/2}2^{-(\log L)^2},
\end{equation}
where we use the inequality $Z_L(\lambda)\geq C_5(\lambda)^{-1}2^L L^{-3/2}$ for all $L\geq 2$ and some $C_5(\lambda)>0$ by Theorem \ref{theorem for the asymptotic behavior of the partition function}.
Secondly,
since  there are at most $L$ corners in any path $\xi\in \Omega_L$, we have
\begin{equation*}
\sum_{x=1}^{L-1} R_x(\xi)\leq L,
\end{equation*} 
where $R_x(\xi)$ is defined in (\ref{the generator of the dynamics}). Therefore, for any subset $\mathcal{A}\subset \Omega_L$ and $s\geq 0$,
\begin{equation}
\mathbb{P}\big[\forall t \in [s, s+L^{-1}]: \sigma_t^{\mu}\in \mathcal{A} \mid \sigma_s^{\mu} \in \mathcal{A}\big]\geq e^{-1}. \label{observation for flippable cornersh}
\end{equation}
Taking $\mathcal{A} \colonequals \{\xi \in \Omega_L: Q(\xi)>(\log L)^2 \}$, we define the occupation time to be
\begin{equation}
u(t) \colonequals \int_0^{t} \mathbbm{1}_{\mathcal{A}}(\sigma_s^{\mu}) ds. \label{occupation time}
\end{equation}
By Fubini's Theorem, we obtain
\begin{equation}
\mathbb{E}[u(2t_{\delta})]=2t_{\delta}\mu(\mathcal{A}). \label{occupation time 1}
\end{equation}
Using (\ref{observation for flippable cornersh}) and strong Markov property, we give a lower bound for $\mathbb{E}[u(2t_{\delta})]$:
\begin{equation}
\mathbb{E}[u(2t_{\delta})] \geq e^{-1}L^{-1} \mathbb{P}[\exists t \in [0, t_{\delta}]: \sigma_t^{\mu}\in \mathcal{A}]. \label{occupation time 2}
\end{equation}
By (\ref{occupation time 1}), (\ref{occupation time 2}) and (\ref{the estimation in the equilibrium state for flippable corners}), we have
\begin{equation}
\mathbb{P}\big[\exists t \in [0, t_{\delta}]: \sigma_t^{\mu}\in \mathcal{A}\big]\leq 2 eLt_{\delta} \mu(\mathcal{A})\leq  4eC_5(\lambda)L^{7/2}t_{\delta}2^{-(\log L)^2}, \label{upper bound using occupation time}
\end{equation}
which vanishes as $L$ tends to infinity. Therefore, we conclude the proof.
\end{proof}

The last ingredient for the proof of Proposition \ref{lemma for controlling the increment for consecutive stopping times} is to control $H(t)$, defined in (\ref{definitionofH}). 
Recall that $t_{\delta}=(1+\delta)\frac{1}{\pi^2}L^2\log L$.
\begin{lemma} \label{Lemma for controlling the highest point}
We have
\begin{equation}
 \lim_{L\rightarrow \infty} \sup_{t\in [t_{\delta/2}, t_{\delta}]}\mathbb{P}\Big[H(t)\geq 2L^{\frac{1}{2}}(\log L)^2\Big]=0.
 \end{equation}
 \end{lemma}
Intuitively, for $\lambda \in (0, 2)$, $\big(\frac{\xi_{[xL]}}{\sqrt{L}} \big)_{x\in [0, 1]}$ under $\muL$ converges to Brownian excursion. Therefore, the dynamics $(\sigma_t^{\wedge})_{t\geq 0}$ is like the simple exclusion process, and we can apply \cite[Theorem 2.4]{lacoin2016mixing} to obtain Lemma \ref{Lemma for controlling the highest point}.  
We postpone the proof in Appendix \ref{proof of the lemma:highestpoint}.
Now, we are ready to prove  Proposition \ref{lemma for controlling the increment for consecutive stopping times}.
\begin{proof}[Proof of Proposition \ref{lemma for controlling the increment for consecutive stopping times}]
 We define the event  $\mathcal{H}_L$ where the highest point of $\sigma_t^{\wedge}$ is not too high and there are a lot of flippable corners in $\sigma_t^{\mu}$ during the time interval $[t_{\delta/2}, t_{\delta/2}+L^2]$, 
\begin{equation*}
\mathcal{H}_L=\bigg\{ \int_{t_{\delta/2}}^{t_{\delta/2}+L^2}   \mathbbm{1}_{\{H(t)\leq 
2L^{\frac{1}{2}}(\log L)^2
\}\bigcap \{Q(\sigma_t^{\mu}) \leq (\log L)^2\}}  dt  \geq L^2 \Big(1-2^{-(K+1)} \Big)   \bigg \}.
\end{equation*}
First, we show that $\mathcal{H}_L$ holds with high probability. We have
\begin{align}
\mathbb{P}\Big[\mathcal{H}^{\complement}_L\Big]&=\mathbb{P}\bigg[\int_{t_{\delta/2}}^{t_{\delta/2}+L^2} \mathbbm{1}_{\{H(t)>2L^{\frac{1}{2}}(\log L)^2\} \bigcup \{Q(\sigma^{\mu}_t)> (\log L)^2\}} dt \geq L^2 2^{-(K+1)} \bigg] \nonumber \\
&\leq \mathbb{P}\bigg[\int_{t_{\delta/2}}^{t_{\delta/2}+L^2} \mathbbm{1}_{\{H(t)>2L^{\frac{1}{2}}(\log L)^2\}}dt \geq L^2 2^{-(K+2)} \bigg]  \nonumber  \\ 
&+  \mathbb{P}\bigg[\int_{t_{\delta/2}}^{t_{\delta/2}+L^2} \mathbbm{1}_{\{Q(\sigma^{\mu}_t)>(\log L)^2\}}dt \geq L^2 2^{-(K+2)} \bigg],  \label{the event H has almost all the mass}
\end{align}
which vanishes as $ L$ tends to infinity, grounded on Markov's inequality,  Lemma \ref{flippable corners for the equilibrium}, Lemma \ref{Lemma for controlling the highest point} and the fact that $K$ is a constant.

From now on, we assume the event $\mathcal{A}_L \cap \mathcal{H}_L \cap\{\Tau_2=t_{\delta/2}\}$.  Based on (\ref{the event H has almost all the mass}), Lemma \ref{the event AL has almost the full mass} and Lemma \ref{flippable corners for the equilibrium}, we have
\begin{equation*}
\lim_{L \to \infty} \mathbb{P}\Big[\mathcal{A}_L \cap \mathcal{H}_L \cap\{\Tau_2=t_{\delta/2}\}\Big]=1.
\end{equation*}
 By induction, we show that $\Delta \Tau_j=\Tau_j-\Tau_{j-1}\leq 2^{-j}L^2$ for all $j \in \lint 3, K \rint$. We argue by contradiction: let  $i_0$ be the smallest integer satisfying
\begin{equation*}
\Delta \Tau_{i_0} > 2^{-i_0}L^2.
\end{equation*}
We know that 
\begin{equation}
\Delta_{i_0}\langle  A \rangle \geq \int_{\Tau_{i_0-1}}^{\Tau_{i_0-1}+2^{-i_0}L^2} \partial_t \langle A. \rangle \mathbbm{1}_{\{ h(t)\leq 2L^{\frac{1}{2}(\log L)^2}\}\cap \{ Q(\sigma^{\mu}_t)\leq (\log L)^2 \}} dt. \label{combine1}
\end{equation}
According to Lemmas \ref{lemma for the lower bound of the derivative of the second weighted area},   \ref{flippable corners for the equilibrium} and \ref{Lemma for controlling the highest point}, we have a lower bound for $\partial_t\langle A.\rangle$ when the indicator function equals to $1$. 
That bound is
\begin{equation}
\partial_t\langle A.\rangle\geq \frac{\lambda  \delta_{\min}}{3(1+\lambda)}\frac{A_t}{H(t)Q(\sigma_t^{\mu})}\geq \frac{\lambda  \delta_{\min}}{6(1+\lambda)}\frac{A_t}{L^{\frac{1}{2}}(\log L)^4}. \label{combine2}
\end{equation}
Since $\Tau_2=t_{\delta/2}$ and $\Delta \Tau_j=\Tau_j-\Tau_{j-1}\leq 2^{-j}L^2$ for $j<i_0$,   we know that
\begin{equation*}
\Tau_{i_0-1}\leq t_{\delta/2}+L^2\sum_{j=3}^{i_0-1}2^{-j}\leq t_{\delta/2}+(1-2^{-(i_0-1)})L^2,
\end{equation*}
and then $\Tau_{i_0-1}+2^{-i_0}L^2\leq t_{\delta/2}+L^2$.
Moreover, when  the assumption $\mathcal{H}_L$ holds, the indicator function 
\begin{equation*}
\mathbbm{1}_{\{ H(t)\leq 2L^{\frac{1}{2}}(\log L)^2\} \cap \{Q(\sigma_t^{\mu})\leq (\log L)^2 \} }
\end{equation*}
is equal to $1$ on a set, which is of Lebesgue measure at least
\begin{equation}
 (2^{-i_0}-2^{-(K+1)})L^2\geq 2^{-(K+1)}L^2. \label{combine3}
\end{equation}
 
Combining (\ref{combine1}), (\ref{combine2}) and (\ref{combine3}), we obtain
\begin{align}
\Delta_{i_0} \langle  A \rangle & \geq 2^{-(K+1)}L^2\frac{\lambda  \delta_{\min}}{6(1+\lambda)}\frac{A_t}{L^{\frac{1}{2}}(\log L)^4}\geq  2^{-(K+1)}\frac{\lambda  \delta_{\min}}{6(1+\lambda)}L^{3-i_0\eta}(\log L)^{-4},  \label{contradiction1}
\end{align}
where the last inequality uses the fact that $A_t>L^{\frac{3}{2}-i_0\eta}$,  for $t<\Tau_{i_0}$.
In addition, since we are in $\mathcal{A}_L$, we know that 
\begin{equation}
\Delta_{i_0}\langle A\rangle \leq L^{3-2(i_0-1)\eta+\frac{1}{2}\eta}. \label{contradiction2}
\end{equation}
However, as $i_0\geq 3$, we have
$$3-2(i_0-1)\eta+\frac{1}{2}\eta<3-i_0\eta.$$ Therefore, there is a contradiction between (\ref{contradiction1}) and (\ref{contradiction2}), as long as $L$ is large enough.
\end{proof}
\section{Upper bound on the mixing time of the dynamics starting from the extremal paths for $\lambda \in (1, 2)$}\label{upper bound for the mixing time of the dynamics starting with extremal paths}
 For the pinning model without a wall (see \cite[Section 1]{caputo2008approach}), the critical value $\lambda_c=1$, while the critical value $\lambda_c=2$ for the pinning model with a wall. Due to the repulsion effect of the wall,
the process $(A_t)_{t\geq 0}$ defined in Subsection \ref{diffusion estimate after contraction} is not a surpermartingale for $\lambda \in (1, 2)$.
But there is still monotonicity in the dynamics starting with the maximal (or minimal) path for $\lambda \in (1, 2)$, which can be exploited to provide an upper bound on the mixing time by
 applying the censoring inequality in \cite[Theorem 1.1]{peres2013can}. This inequality says that canceling some prescribed updates slows down the mixing of the Glauber dynamics starting from the maximal (or minimal) configuration of a monotone spin system. 
 
 Let us state the setting for applying the censoring inequality.
 A censoring scheme is a c\`adl\`ag function defined by
\begin{equation*}
\mathcal{C} \colon \mathbb{R}^+ \to \mathcal{P}(\Theta),
\end{equation*}
where $\Theta$ is defined in (\ref{the set of spins for censoring}) and $\mathcal{P}(\Theta)$ is the set of all subsets of $\Theta$.
The censored dynamics with a censoring scheme $\mathcal{C}$ is the dynamics obtained from the graphical construction in Subsection \ref{graphical construction}, except that the update  at time $t$ is canceled if and only if it is an element of $\mathcal{C}(t)$.  In other words,  we construct the dynamics by using the graphical construction in Subsection \ref{graphical construction} with one extra rule: 
if $\Tau_{(x,z)}^{\uparrow}$ or $\Tau_{(x,z)}^{\downarrow}$ rings at time $t$, the update is performed if and only if $(x,z) \not\in \mathcal{C}(t)$. Let  $(\sigma_t^{\xi, \mathcal{C}})_{t\geq 0}$ denote the trajectory of the censored  dynamics  with a censoring scheme $\mathcal{C}$ and starting from the path $\xi \in \Omega_L$, and let $P_t^{\xi, \mathcal{C}}$ denote  the law of distribution of the time marginal $\sigma_t^{\xi, \mathcal{C}}$.

 The Glauber dynamics of this polymer pinning model is a monotone spin system in the sense of \cite[Subsection 1.1]{peres2013can} (detailed in Appendix \ref{spinsystem}), and we refer to Figure \ref{spinsystemfig} in Appendix \ref{spinsystem}  for a quick look.
The following proposition follows directly from \cite[Theorem 1.1]{peres2013can}.
\begin{proposition}\label{Peres-Winkler inequality}
For any prescribed censoring scheme $\mathcal{C}$, for all $\lambda\in [0, \infty)$, all $t\geq 0$ and $\xi \in \big\{\wedge,\vee \big\}$, we have
\begin{equation}
\Vert P_t^{\xi}-\mu \Vert_{\TV} \leq \Vert P_t^{\xi, \mathcal{C}}-\mu\Vert_{\TV}.
\end{equation}
\end{proposition}

Besides Proposition \ref{Peres-Winkler inequality}, we need the two following results in the proof of the upper bound on the mixing time.
Firstly, by \cite[Lemmas 20.5 and 20.11]{LPWMCMT}, we know  that the asymptotic rate of convergence to equilibrium of this reversible Markov chain is 
\begin{equation}
\lim_{t \rightarrow \infty} t^{-1} \log d^{L, \lambda}(t)=-\gap_L,  \label{definition for the spectral gap}
\end{equation}
where $\gap_L>0$ is the smallest nonzero eigenvalue of $-\mathcal{L}$, usually referred to as the spectral gap.
By monotonicity of the Glauber dynamics and (\ref{upper bound for mixing time in terms of coupling time}), for all $\lambda>0$ we have
\begin{equation}
d^{L, \lambda}(t)\leq \mathbb{P}\Big(\sigma_t^{\wedge}\neq \sigma_t^{\vee}\Big)=\mathbb{P}\Big(\Phi(\sigma_t^{\wedge})-\Phi(\sigma_t^{\vee})\geq 2\sin(\tfrac{\pi}{L})\Big),\label{inequality for total variation distance in terms of monotone function and minimal increment}
\end{equation}
where $\Phi(\xi)$ is defined in (\ref{weighted area function introduced by Caputo}). Moreover, for all $\lambda>0$, by \cite[Equation (4.1)]{caputo2008approach}  we have
 $$ \mathbb{E}\big[\Phi(\sigma_t^{\wedge})\big]-\mathbb{E}\big[\Phi(\sigma_t^{\vee})\big]\leq \tfrac{L^2}{2}e^{-t \kappa_L}.$$
Applying Markov's inequality, we reclaim the useful result in \cite{caputo2008approach}.  
\begin{lemma}For all $\lambda>0$, we have
\begin{equation}
d^{L, \lambda}(t) \leq \frac{L^2 e^{-\kappa_L t}}{4\sin(\frac{\pi}{L})}. \label{an upper bound for the mixing time in Caputo}
\end{equation}
Plugging this into (\ref{definition for the spectral gap}), we obtain
 \begin{equation}
 \gap_L\geq \kappa_L= 1-\cos \Big(\frac{\pi}{L}\Big).  \label{lower bound for the spectral gap}
 \end{equation}
\end{lemma}

Secondly, the following lemma is an application of the Cauchy-Schwarz inequality and the reversibility of the Markov chain. For reference, we mention \cite[Equation (2.6)]{caputo2012polymer}. 
\begin{lemma}\label{the standard inequality for the total variation distance by using L2 bound and the spectral gap} 
For any probability distribution $\nu$ on $\Omega_L$,  we have
\begin{equation}
\Vert \nu P_t -\mu \Vert_{\TV}\leq \frac{1}{2} e^{-t \cdot \gap_L }\sqrt{\Var_{\mu}(\rho)},
\end{equation}
where $\rho\colonequals\frac{d \nu}{d \mu}$ and $\Var_{\mu}(\rho)\colonequals\mu(\rho^2)-\mu(\rho)^2$.
\end{lemma}

We define 
\begin{equation}
G_L \colonequals \big\{(x, 1): x \in \lint 2, L-2 \rint \cap 2 \mathbb{N} \big\} \label{the green squares}
\end{equation}
where $\Tau_{(x,1)}^{\uparrow}$ or $\Tau_{(x,1)}^{\downarrow}$ rings, the update\textemdash in the graphical construction of Subsection \ref{graphical construction}\textemdash changes the number of contact points $\mathcal{N}$, defined in (\ref{number of contact points}). Moreover, $G_L$ corresponds to the centers of the green squares shown in Figure \ref{spinsystemfig}.
Before we start the proof of the upper bound on the mixing time for the dynamics starting with the maximal path $\wedge$,  we outline the idea.
\begin{itemize}
\item[(i)] We elaborate a censoring scheme $\mathcal{C}$, where $\mathcal{C}(t)=G_L$ for $t< t_{\delta/2}$ and $\mathcal{C}(t)= \varnothing$ 
for $t\geq t_{\delta/2}$.
We run the dynamics $(\sigma_t^{\wedge, \mathcal{C}})_{0\leq t<t_{\delta/2}}$ with this censoring scheme.
\item[(ii)] By Theorem \ref{the cutoff phenomenon in the diffusive regime}, the distribution of $\sigma_{t_{\delta/2}}^{\wedge,\mathcal{C}}$ is close to $\mu_{L}^{0}$ in total variation distance. 
\item[(iii)] As the Radon-Nikodym derivative of  $\mu_L^0$ with respect to $\mu_{L}^{\lambda}$ is bounded by a constant, we apply Lemma \ref{the standard inequality for the total variation distance by using L2 bound and the spectral gap} and use (\ref{lower bound for the spectral gap}) to conclude the proof.
\end{itemize}

\begin{proposition}\label{proposition for the mixing time starting with the maximal path in the second regime}
For any $\lambda \in (1, 2)$, any $\epsilon>0$ and any $\delta>0$,  if $L$ is sufficiently large, we have 
\begin{equation}
\Tmaximal(\epsilon)\leq \frac{1+\delta}{\pi^2}L^2 \log L.
\end{equation}
\end{proposition}

\begin{proof}
Recall that $\mathcal{N}$ is the number of contact points, defined in (\ref{number of contact points}).
We run the dynamics starting from the maximal path $\wedge$,  censoring those updates which change the value of contact points $\mathcal{N}$ for $t<t_{\delta/2}$. More precisely,
recalling $t_{\delta}=(1+\delta)\frac{1}{\pi^2}L^2 \log L$, we present a censoring scheme $\mathcal{C}\colon \mathbb{R}^+ \to \mathcal{P}(\Theta)$, defined by
\begin{equation*}
\mathcal{C}(t)\colonequals
\begin{cases*} G_L  & if  $t \in [0, t_{\delta/2})$,\\
 \varnothing  & if $t \in [t_{\delta/2}, \infty)$.
\end{cases*}
\end{equation*}
We recall that $\sigma_t^{\wedge, \mathcal{C}}$ is the dynamics constructed by using the graphical construction with one extra rule: when the clock process $\Tau_{(x,1)}^{\uparrow}$ or $\Tau_{(x, 1)}^{\downarrow}$ rings for any $x\in \lint 1, L-1\rint \cap 2\mathbb{N}$ and all $t<t_{\delta/2}$, we do not update. We refer to Figure \ref{fig:upperboundmaximalpath} for illustration. While $t\geq t_{\delta/2}$, $(\sigma_t^{\wedge, \mathcal{C}})_{t\geq t_{\delta/2}}$ is constructed by the graphical construction in Subsection \ref{graphical construction} without censoring.

\begin{figure}[h]
\centering
  \begin{tikzpicture}[scale=.4,font=\tiny]
   \draw (25,4) -- (25,-1) -- (52,-1);
    \draw[color=blue] (25,-1)--(26,0) -- (27,1) -- (28,0) -- (29,1) -- (30,2) -- (31,3) -- (32,4) -- (33,3) -- (34,2) -- (35,1) -- (36,0) -- (37,1) -- (38,2) -- (39,1) -- (40,0) -- (41,1) -- (42,0) -- (43,1) -- (44,2) -- (45,1) -- (46,0) -- (47,1) -- (48,2) -- (49,1) -- (50,0)--(51, -1);
    \foreach \x in {25,...,51} {\draw (\x,-1.3) -- (\x,-1);}
    \draw[color=red,dashed] (40,0) -- (41,-1) -- (42,0); \draw (41,0.8) edge[out=290, in=70, ->] (41,-0.8);\node[color=red,below] at (41.4, 0.4) {$\times$};
    \draw[color=red,dashed](26,0)--(27,-1)--(28,0); \draw (27,0.8) edge[out=290, in=70, ->] (27,-0.8); \node[color=red,below] at (27.4, 0.4) {$\times$};
   \draw[dashed] (35,1) -- (36,2) -- (37,1); \draw (36,0.2) edge[out=70, in=290, ->] (36,1.8); \node[below] at (36.4, 1.6) {$\frac{1}{2}$};
    \draw[dashed] (31,3) -- (32,2) -- (33,3);\draw (32,3.8) edge[out=290, in=70, ->] (32,2.2); \node[below] at (32.4, 3.6) {$\frac{1}{2}$};
   \draw[dashed](43,1)--(44,0)--(45,1);\draw (44,1.8) edge[out=290, in=70, ->] (44,0.2);\node[below] at (44.4, 1.6) {$\frac{1}{2}$};
    \draw[fill] (25,-1) circle [radius=0.1];
    \draw[fill] (51,-1) circle [radius=0.1];   
    \node[below] at (25,-1.3) {$0$};
    \node[below] at (51,-1.3) {$L$};
    \node[blue,above] at (38, 1.8){$\sigma_t^{\wedge,\mathcal{C}}$};
    \draw[thick,->] (25,-1) -- (25,4) node[anchor=north west]{y};
    \draw[thick,->] (25,-1) -- (52,-1) node[anchor=north west]{x};
  \end{tikzpicture}
  \caption{\label{fig:upperboundmaximalpath}  A graphical representation of the jump rates for the dynamics $\sigma_t^{\wedge, \mathcal{C}}$ when $t<t_{\delta/2}$. Those red dashed corners are not available and  labeled with $\times$, while the other corners are flippable with rate $1/2$.  }
\end{figure}
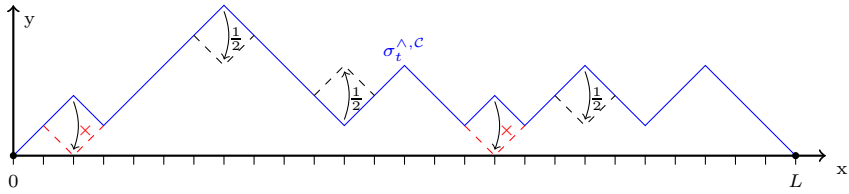

Now we show that $P_{t_{\delta/2}}^{\wedge, \mathcal{C}}$ is close to $\mu_{L}^0$.
By Remark \ref{remark}, applying Theorem \ref{the cutoff phenomenon in the diffusive regime},
for all $\lambda \in (1, 2)$,  all $\delta>0$ and all $\epsilon>0$,  if $L$ is sufficiently large,  we have 
\begin{equation}
\Vert P_{t_{\delta/2}}^{\wedge, \mathcal{C}} -\mu_L^0 \Vert_{\TV} \leq \epsilon/2. \label{the dynamics of the maximal path with the cencoring scheme is already well mixed in the space without contact point}
\end{equation}
For any $\xi \in \Omega_L$, define
$$\rho (\xi)\colonequals\frac{ d\mu_L^0}{d \muL}(\xi),$$ and we want to show that $\rho$ is bounded above uniformly for $\xi \in \Omega_L$. 
For any $\xi\in \Omega_L\setminus \Omega_L^{+}$\textemdash recalling $\Omega_L^+=\{\xi \in \Omega_L: \mathcal{N}(\xi)=0 \}$, since $\mu_L^0(\xi)=0$,  $$\rho(\xi)= \frac{ \mu_L^0(\xi)}{ \muL(\xi)}=0.$$
While for any $\xi \in \Omega_L^{+}$, applying Theorem \ref{theorem for the asymptotic behavior of the partition function},  for all $L\geq 4$ we have
\begin{equation*}
\rho(\xi)=\frac{d \mu_L^0}{d \muL}(\xi)=\frac{\mu_L^0(\xi)}{\muL(\xi)}=\frac{1/Z_{L-2}(1)}{1/Z_L(\lambda)}\leq C_5(\lambda),  
\end{equation*}
where $C_5(\lambda)>0$ is a suitable constant and only depends on $\lambda$.
By Lemma \ref{the standard inequality for the total variation distance by using L2 bound and the spectral gap} and 
 (\ref{lower bound for the spectral gap}), for any given $\delta>0$, we have
 \begin{equation}
 \lim_{L\to \infty} \Vert \mu_L^0P_{\frac{\delta}{2} L^2\log L}-\muL \Vert_{\TV} =0.  \label{upperboundfreedom} 
 \end{equation}
 
At this moment,  we are ready to show that $P_{t_{\delta}}^{\wedge, \mathcal{C}}$\textemdash the  distribution of the censored dynamics at $t_{\delta}$\textemdash is close to the stationary measure $\muL$.
 By the definition of $\mathcal{C}$, we have
\begin{align*}
\Vert P_{t_{\delta}}^{\wedge, \mathcal{C}} -\muL \Vert_{\TV}&=\Vert P_{t_{\delta/2}}^{\wedge, \mathcal{C}}P_{t_{\delta}-t_{\delta/2}} -\muL \Vert_{\TV}\\
&\leq \Vert P_{t_{\delta/2}}^{\wedge, \mathcal{C}} P_{t_{\delta}-t_{\delta/2}}-\mu_L^0P_{t_{\delta}-t_{\delta/2}}\Vert_{\TV}+\Vert \mu_L^0P_{t_{\delta}-t_{\delta/2}}-\muL\Vert_{\TV} \\
&\leq  \Vert P_{t_{\delta/2}}^{\wedge, \mathcal{C}} -\mu_L^0\Vert_{\TV}+\Vert \mu_L^0P_{t_{\delta}-t_{\delta/2}}-\muL\Vert_{\TV}.  \numberthis \label{upper bound for the total variation distance for the maximal path}
\end{align*}
Here the first inequality uses the triangle inequality.  The second inequality is based on the fact that $\Vert \alpha P_t -\beta P_t\Vert_{\TV}\leq \Vert \alpha -\beta \Vert_{\TV}$ for any two probability measures $\alpha, \beta$ on $\Omega_L$, and $P_t$ is a transition matrix on $\Omega_L$. 
The first term in (\ref{upper bound for the total variation distance for the maximal path}) is not bigger than $\epsilon/2$ by (\ref{the dynamics of the maximal path with the cencoring scheme is already well mixed in the space without contact point}) for $L$ sufficiently large.
The second term in (\ref{upper bound for the total variation distance for the maximal path}) is smaller than or equal to $\epsilon/2$ by (\ref{upperboundfreedom}) for $L$ sufficiently large.

 Recall that $P_t^{\wedge}$ is the distribution of $\sigma_t^{\wedge}$ without censoring. By Proposition \ref{Peres-Winkler inequality},  for any $t\geq 0$, we have
\begin{equation}
\Vert P_t^{\wedge} -\muL \Vert_{\TV} \leq \Vert P_t^{\wedge, \mathcal{C}} -\muL \Vert_{\TV}. \label{application of the censoring inequality}
\end{equation}
 Combining (\ref{upper bound for the total variation distance for the maximal path}) and (\ref{application of the censoring inequality}), we conclude the proof.
\end{proof}

Our next task is to provide an upper bound on the mixing time for the dynamics starting from the minimal path.
\begin{proposition}\label{upper bound for the minimal path in the second diffusive regime}
For any $\lambda \in (1, 2)$, any $\epsilon>0$ and any $\delta>0$, if  $L$ is sufficiently large,  we have
\begin{equation}
\Tminimal(\epsilon)\leq \frac{1+\delta}{\pi^2}L^2\log L.
\end{equation}
\end{proposition}
The idea for the proof of Proposition \ref{upper bound for the minimal path in the second diffusive regime} is similar to Proposition \ref{proposition for the mixing time starting with the maximal path in the second regime}.
We first prove that
  by time $s_0(L) \colonequals 10L^{16/9} \log L$, the contact points between the $x-$axis and the dynamics starting from the minimal path are close to the boundaries, which is the following lemma.

\begin{lemma}\label{lemma for only contact points close to the boundary for the dynamics starting from the minimal path}
For any given $\epsilon>0$ and $\lambda \in (1, 2)$,  let $M=M(\lambda, \epsilon)$ be a positive integer, and
\begin{equation}
\mathcal{E}_{L,M} \colonequals \Big\{ \xi \in \Omega_L: \xi_x \geq 1, \forall x \in \lint M, L-M\rint  \Big\}. \label{definition of paths with contact points close to boundaries}
\end{equation}
For all $L\geq 2M$,  we have
\begin{equation}
\mathbb{P}\Big[\sigma_{s_0}^{\vee} \in \mathcal{E}_{L,M}\Big]\geq 1-\epsilon/2.
\end{equation}
\end{lemma}

\begin{proof}
Let $m$ and $n$ be two positive integers, and $n<m<L/2$.  Observe that in the graphical construction,  if we run  the dynamics $(\sigma_t^{\vee})_{t\geq 0}$ with  the points $(2n, 0)$ and $(2m, 0)$ fixed for all $t\geq 0$, denoted as $(\overline{\sigma}_t^{\vee})_{t\geq 0}$ with $\overline{\sigma}_t^{\vee}(2n)\equiv \overline{\sigma}_t^{\vee}(2m)\equiv 0$ for all $t\geq 0$,  we have
\begin{equation}
 \forall t\geq 0, \mbox{ }\overline{\sigma}_t^{\vee}\leq \sigma_t^{\vee}. \label{stochatic domination}
\end{equation}
By symmetry, to give an upper bound on $\mathbb{P}[\sigma_t^{\vee}(x)=0]$, we only need to consider $x\in \lint 0, L/2\rint$.
For all $M\leq x \leq L/2$ and $x\in 2\mathbb{N}$, let $\bar{x}\colonequals 2\lfloor x^{8/9}/2 \rfloor$ and $\bar{L} \colonequals 2\lfloor L^{8/9}/2\rfloor$. 
Let $d^{L, \lambda}(t)\leq L^{-3/2}$ in (\ref{an upper bound for the mixing time in Caputo}). For all $L$ sufficiently large, we obtain
\begin{equation}
\T(L^{-3/2})\leq \frac{18}{\pi^2}L^2\log L .   \label{upper bound for the mixing time for choose epsilon sufficiently small depending on the system}
\end{equation} 
 Therefore, the quantity $s_0$ satisfies
\begin{equation*}
T_{\mathrm{mix}}^{\bar{L}} \bigg( \bar{L}^{-3/2}\bigg)\leq s_0.
\end{equation*}
 Using (\ref{stochatic domination}), (\ref{upper bound for the mixing time for choose epsilon sufficiently small depending on the system}) and (\ref{estimation for a hitting of the origin}) respectively, for all $t\geq s_0$,  we obtain
\begin{align*}
\mathbb{P}\Big[\sigma_t^{\vee}(x)=0\Big]&\leq \mathbb{P}\Big[\sigma_t^{\vee}(x)=0 \vert  \sigma_t^{\vee}(x-\bar{x}) \equiv \sigma_t^{\vee}(x+\bar{x})\equiv 0, \forall  t\geq 0 \Big]\\
&\leq \mu_{2\bar{x}}^{\lambda}(\xi_{\bar{x}}=0)+\Vert P_t^{\vee} -\mu_{2\bar{x}}^{ \lambda} \Vert_{\TV}\\
& \leq C_6(\lambda)x^{-4/3},\numberthis \label{upper bound for the dynamics hitting the origin after a sufficiently big time}
\end{align*}
where $C_6(\lambda)>0$ only depends on $\lambda$.
In the second inequality, there is an abuse of notation\textemdash $P_t^{\vee}$ denotes the distribution of $\sigma_t^{\vee}$ starting with  the minimal path $\vee$ of $\Omega_{2\bar{x}}$.  
Therefore, due to symmetry and (\ref{upper bound for the dynamics hitting the origin after a sufficiently big time}), we obtain
\begin{equation}
\begin{aligned}
\sum_{x=M}^{L-M} \mathbb{P}[\sigma_{s_0}^{\vee}(x)=0]&=2\sum_{x=M}^{L/2} \mathbb{P}[\sigma_{s_0}^{\vee}(x)=0] \\&
\leq 2C_7(\lambda)M^{-1/3}. \label{the choice of the constant M}
\end{aligned}
\end{equation}
Let $C(\lambda, \epsilon)>0$ be a constant such that the right-hand side is smaller than $\epsilon/2$, if $M\geq C(\lambda, \epsilon)$.
Applying Markov's inequality and (\ref{the choice of the constant M}), we obtain
\begin{equation}
\mathbb{P}\Big[\sigma_{s_0}^{\vee} \notin \mathcal{E}_{L,M}\Big]= \mathbb{P}\Big[\sum_{x=M}^{L-M} \mathbbm{1}_{ \{\sigma_{s_0}^{\vee}(x)=0\}}\geq 1\Big]\leq \epsilon/2  . \label{the event B_L has almost all the mass}
\end{equation}
\end{proof}

Now we show that the dynamics starts from $\xi\in \mathcal{E}_{L, M}$, and censors the updates that change the number of the contact points until time $t_{\delta/2}$. Its distribution is close to $\muL$ in total variation distance at time $t_{3\delta/4}$.

\begin{lemma}\label{mixing time for brownian excursion type paths}
 Let $\xi \in \mathcal{E}_{L, M}$, and let $(\sigma_t^{\xi, \mathcal{C}})_{t\geq 0}$ be a censored dynamics with the censoring scheme $\mathcal{C}\colon \mathbb{R}^+ \to \mathcal{P}(\Theta)$ defined by
\begin{equation*}
\mathcal{C}(t)\colonequals \begin{cases*}  
G_L & if $ t\in [0, t_{\delta/2})$,\\
 \varnothing & if $t\in [t_{\delta/2}, \infty)$.
\end{cases*}
\end{equation*}
where $G_L$ is defined in (\ref{the green squares}). For any given $\epsilon>0$, for all $L$ sufficiently large, we have
\begin{equation}
\Vert P_{t_{3\delta/4}}^{\xi, \mathcal{C}}-\muL\Vert<\epsilon/2,
\end{equation}
where we recall that $t_{\delta}=(1+\delta)\pi^{-2}L^2 \log L$ and $P_t^{\xi, \mathcal{C}}$ denotes the marginal distribution of the censored dynamics $(\sigma_t^{\xi, \mathcal{C}})_{t\geq 0}$ at time t.
\end{lemma}

With Lemma \ref{mixing time for brownian excursion type paths} at hand, we are ready to prove Proposition \ref{upper bound for the minimal path in the second diffusive regime}.
Combining  Lemma \ref{lemma for only contact points close to the boundary for the dynamics starting from the minimal path}, Lemma \ref{mixing time for brownian excursion type paths} and Proposition \ref{Peres-Winkler inequality}, we conclude the proof of Proposition \ref{upper bound for the minimal path in the second diffusive regime}, since
$s_0+t_{3\delta/4}\leq t_{\delta}$.

\begin{proof}[Proof of Lemma \ref{mixing time for brownian excursion type paths}]
 For $\xi \in \mathcal{E}_{L,M}$, set
\begin{equation}
\begin{aligned}
&l(\xi)\colonequals\sup \big\{x\leq M: \xi_x=0  \big\},\\
&r(\xi)\colonequals\inf \big\{ x\geq L-M:  \xi_x=0 \big\}.
\end{aligned} \label{observation of independence of the three intervals}
\end{equation}
Observe that  the censored dynamics $(\sigma_t^{\xi, \mathcal{C}})_{0\leq t<t_{\delta/2}}$ restricted in the intervals $\lint 0, l  \rint$, $\lint l, r\rint$ and $\lint r, L \rint$ respectively  are independent. Let the marginal distribution restricted in these three intervals be denoted by $P_{t, l}^{\xi, \mathcal{C}}$, $P_{t, r-l}^{\xi, \mathcal{C}}$, $P_{t, L-r}^{\xi, \mathcal{C}}$ respectively. 
 We refer to Figure \ref{fig:forminimalpath} for illustration. 
\begin{figure}[h]
\centering
  \begin{tikzpicture}[scale=.4,font=\tiny]
   \draw (25,4) -- (25,-1) -- (52,-1);
    \draw[color=blue] (25,-1)--(26,0) -- (27,1) -- (28,2) -- (29,1) -- (30,0) -- (31,-1) -- (32,0) -- (33,1) -- (34,0) -- (35,1) -- (36,2) -- (37,1) -- (38,0) -- (39,1) -- (40,0) -- (41,1) -- (42,0) -- (43,1) -- (44,0) -- (45,1) -- (46,0) -- (47,-1) -- (48,0) -- (49,-1) -- (50,0)--(51, -1);
    \foreach \x in {25,...,51} {\draw (\x,-1.3) -- (\x,-1);}
    \draw[fill] (25,-1) circle [radius=0.1];
    \draw[fill] (51,-1) circle [radius=0.1];   
    \node[below] at (25,-1.3) {$0$};
    \node[below] at (51,-1.3) {$L$};
    \draw[fill, color=green] (32,-1) circle [radius=0.1];
   \node[below]at (32,-1.3){$M$};
   \draw[fill, color=green] (44,-1) circle [radius=0.1];
     \node[below]at (44,-1.3){$L-M$};
    \draw[fill, color=red] (31,-1) circle [radius=0.1];
     \node[below]at (31,-1.3){$l$};  
      \draw[fill, color=red] (47,-1) circle [radius=0.1];
     \node[below]at (47,-1.3){$r$};  
    \draw[thick,->] (25,-1) -- (25,4) node[anchor=north west]{y};
    \draw[thick,->] (25,-1) -- (52,-1) node[anchor=north west]{x};
    \draw[dashed] (27,1) -- (28,0) -- (29,1); \draw (28,1.8) edge[out=290, in=70, ->] (28,.2);     \node[below] at (28.5, 1.8) {$\tfrac{1}{2}$};
    \draw[dashed](35,1) -- (36,0) -- (37,1); \draw (36,1.8) edge[out=290, in=70, ->] (36,.2);     \node[below] at (36.5, 1.8) {$\tfrac{1}{2}$};
    \draw[color=red, dashed](30,0) -- (31,1) -- (32,0); \draw (31,-0.8) edge[out=70, in=290, ->] (31,.8);     \node[color=red,below] at (31.4, .7) {$\times$};
    \draw[color=red, dashed](46,0) -- (47,1) -- (48,0); \draw (47,-0.8) edge[out=70, in=290, ->] (47,.8);     \node[color=red,below] at (47.4, .7) {$\times$};
      \draw[color=red, dashed](48,0) -- (49,1) -- (50,0); \draw (49,-0.8) edge[out=70, in=290, ->] (49,.8);     \node[color=red,below] at (49.4, .7) {$\times$};
    \draw[dashed, color=red](40,0) -- (41,-1) -- (42,0); \draw (41,0.8) edge[out=290, in=70, ->] (41,-.8);     \node[color=red,below] at (41.4, .7) {$\times$};
    \draw[dashed] (43,1) -- (44,2) -- (45,1);\draw (44,0.2) edge[out=70, in=290, ->] (44,1.8); \node[below] at (44.4, 1.8) {$\tfrac{1}{2}$};
   \node[above, color=blue] at (42, 1.2){$\sigma_t^{\xi, \mathcal{C}}$};
  \end{tikzpicture}
  \caption{\label{fig:forminimalpath}  A graphical representation of the jump rates for the censored dynamics $(\sigma_t^{\xi, \mathcal{C}})_{0\leq t< t_{\delta/2}}$ starting from $\xi \in \mathcal{E}_L$. The red dashed corners are not available corners, labeled with $\times$.  To the left hand side of the green point $(M, 0)$, the red point $(l, 0)$ is the first contact point with the $x$-axis at time $t=0$. Moreover, the corner at $(l, 0)$ is fixed for $t\in [0, t_{\delta/2})$. Likewise, the same phenomenon holds for the green point $(L-M, 0)$ and the red point $(r,0)$. In the time interval $[0, t_{\delta/2})$, the censored dynamics $(\sigma_t^{\xi, \mathcal{C}})_{0\leq t< t_{\delta/2}}$ does not touch the $x-$axis in the interval $\lint l+1, r-1\rint$.  
 }
\end{figure}
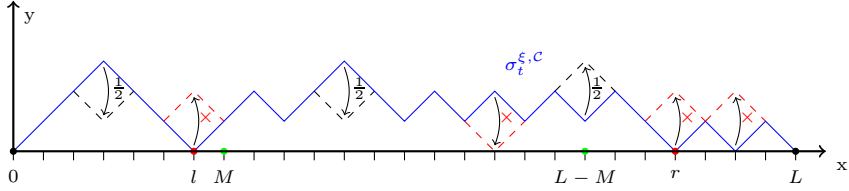

Let the censored  dynamics restricted in the interval $\lint l, r \rint$ be denoted by $(\widetilde{\sigma}^{\xi}_t)_{t<t_{\delta/2}}$, whose invariant probability measure is $\mu_{r-l}^{0}$, defined in (\ref{the stationary measure}). 
By Theorem \ref{the cutoff phenomenon in the diffusive regime} and Remark \ref{remark}, for given $\delta>0$ and $\epsilon>0$, for all $L$ sufficiently large, we have
\begin{equation}
\Vert P_{t_{\delta/2}, r-l}^{\xi, \mathcal{C}}- \mu_{r-l}^0 \Vert_{\TV}\leq \epsilon/4. \label{upper bound for the total variation distance for the excursion path}
\end{equation}
Note that the upper bound in (\ref{upper bound for the total variation distance for the excursion path}) does not depend on the value of $(l, r)$.
Moreover, observe that for any $\xi'\in \Omega_L$,
the product distribution $P_{t_{\delta/2},l}^{\xi, \mathcal{C}}\otimes \mu_{r-l}^0 \otimes P_{t_{\delta/2},L-r}^{\xi, \mathcal{C}}$ satisfies 
\begin{equation}
\Big( P_{t_{\delta/2},l}^{\xi, \mathcal{C}}\otimes \mu_{r-l}^0 \otimes P_{t_{\delta/2},L-r}^{\xi, \mathcal{C}}\Big) (\xi')\leq  \frac{1}{Z_{r-l}(0)},
\end{equation}
while $\mu_L^{\lambda}(\xi')\geq 1/Z_L(\lambda)$ since $\lambda \in (1, 2)$. Therefore, for all $L>2M$ and for any $\xi' \in \Omega_L$, we have
\begin{equation}
\frac{d P_{t_{\delta/2},l}^{\xi, \mathcal{C}}\otimes \mu_{r-l}^0 \otimes P_{t_{\delta/2},L-r}^{\xi, \mathcal{C}}} {d \mu_L^{\lambda}}(\xi')\leq C_8(\lambda)2^{2M}, \label{preparation for radon-nikodym derivative is bounded }
\end{equation}
where the last inequality uses Theorem \ref{theorem for the asymptotic behavior of the partition function} and $r-l\geq L-2M$, since $\xi \in \mathcal{E}_{L,M}$. Note that the right-most hand side in (\ref{preparation for radon-nikodym derivative is bounded }) does not depend on the value of $(l, r)$. 
 We repeat (\ref{upper bound for the total variation distance for the maximal path}) to obtain that for all $L$ sufficiently large,
\begin{equation}
\Vert P_{t_{3\delta/4}}^{\xi, \mathcal{C}}-\muL \Vert \leq \epsilon/2. \label{total variational distance for the minimal path with given condition}
\end{equation}
\end{proof}
Theorem \ref{the cutoff for the extremal paths} is a combination of Proposition \ref{lower bound of the mixing time}, Proposition \ref{proposition for the mixing time starting with the maximal path in the second regime}, and Proposition \ref{upper bound for the minimal path in the second diffusive regime}.

\begin{appendix}
\section{Proof of lemma \ref{Lemma for controlling the highest point}. }\label{proof of the lemma:highestpoint}
 We lift the maximal path $\wedge$ up by a height $L^{1/2}(\log L)^2$. To be precise, define  $\bwedge \colonequals \wedge+m$, \textit{i.e.} $\bwedge_x=\wedge_x+m$ for all $x\in \lint 0, L\rint$, where  $m \colonequals 2\lceil L^{1/2}(\log L)^2/2 \rceil$. 
The graphical construction in Subsection \ref{graphical construction}, with $\Theta$ changed to be 
$$\Theta' \colonequals \Big\{(x, z):  x\in \lint 1, L-1\rint,  z\in \lint 1, m +L/2-1-\vert x- L/2 \vert \rint, x+z \in 2\mathbb{N}+1 \Big\},$$
 allows us to couple the three dynamics $(\sigma_t^{\wedge,\lambda})_{t\geq 0}$, $(\sigma_t^{\barwedge, \lambda})_{t\geq 0}$ and $(\sigma_t^{\barwedge, 0})_{t\geq 0}$, starting from $\wedge$, $\barwedge$ and $\barwedge$ respectively, with parameter $\lambda$, $\lambda$ and $0$ respectively.
By the monotonicity of the starting paths and the parameters $\lambda$ in the dynamics, asserted in Proposition \ref{preserving the monotonicity}, we have 
\begin{equation*}
\begin{aligned} 
 &\sigma_t^{\wedge, \lambda}\leq  \sigma_t^{\bwedge, \lambda},\\
 & \sigma_t^{\bwedge, \lambda}\leq \sigma_t^{\bwedge, 0}.
\end{aligned}
\end{equation*}
 Set 
 \begin{equation*}
 \barh(t)\colonequals\max_{x \in \lint 0, L \rint}\sigma_t^{\bwedge, 0}(x).
 \end{equation*}
 Since $\barh(t)\geq H(t)$, it is enough to prove that
 \begin{equation}
 \lim_{L\rightarrow \infty} \mathbb{P} \Big[ \exists  t\in [t_{\delta/2}, t_{ \delta}]:\barh(t)\geq 2L^{1/2}(\log L)^2\Big]=0, \label{alternativeupperboundhighestpoint}
 \end{equation}
 where we recall that  $t_{\delta}=(1+\delta) \frac{1}{\pi^2}L^2 \log L$.
 We obtain such an upper bound in (\ref{alternativeupperboundhighestpoint}) by comparing $(\sigma_t^{\bwedge,0})_{t\geq 0}$ with the symmetric simple exclusion process.

\subsection{Simple exclusion process.} Define
\begin{equation}
\mathcal{S}_L\colonequals \Big\{\zeta \in \mathbb{Z}^{L+1}: \zeta_0=\zeta_L=m; \vert \xi_{x+1}-\xi_{x} \vert=1, \forall x\in \lint 0, L-1 \rint \Big\}, \label{the height function space of simple exclusion process}
\end{equation} 
and $$ \mathcal{S}_L^+ \colonequals \Big\{ \zeta \in \mathcal{S}_{L}: \zeta_x\geq 1, \forall x \in \lint 0, L \rint  \Big\}.$$
 We define a Markov chain on 
  $\mathcal{S}_L$ by specifying its generator $\mathfrak{L}$. The generator $\mathfrak{L}$ is defined by its action on the functions $\mathbb{R}^{\mathcal{S}_L}$,
\begin{equation}
(\mathfrak{L}f)(\zeta)\colonequals\frac{1}{2}\sum_{x=1}^{L-1}\Big(f(\zeta^x)-f(\zeta)\Big), \label{generator of the exclusion process}
\end{equation}
where $\zeta^x \in \mathcal{S}_L$ is defined by
\begin{equation*}
\zeta^x_y\colonequals\begin{cases*}
\zeta_y &if $y\neq x$,\\
\zeta_{x-1}+\zeta_{x+1}-\zeta_x & if $y=x$.
\end{cases*}
\end{equation*}
When $\zeta_{x-1}=\zeta_{x+1}$, $\zeta$ displays a local extremum at $x$ and we obtain $\zeta^x$ by flipping the corner of $\xi$ at the coordinate $x$.
Let $U_L$ denote the uniform probability measure on $\mathcal{S}_{L}$.  We can see that this Markov chain is reversible with respect to the uniform measure $U_L$. Therefore, $U_L$ is the invariant probability measure for this Markov chain.
 The Markov chain starting with the maximal path $\bwedge$ is denoted by $(\eta^{\bwedge}_t)_{t\geq 0}$. Likewise, let $(\eta_t^{U_L})_{t\geq 0}$ denote the Markov chain with generator $\mathfrak{L}$ and starting path chosen by sampling $U_L$.
 There is a one-one correspondence between this Markov chain and the symmetric simple exclusion process, for which we refer to \cite[Section 2.3]{lacoin2016mixing} for more information.  
Under the measure $U_L$, typical path $\zeta\in \mathcal{S}_L$ does not touch the $x$-axis, which is the following lemma.

\begin{lemma}\label{upper bound for the lifted up system to hit the zero level line}
 For all $L$ sufficiently large, we have
\begin{equation}
U_L(\mathcal{S}_{L}\setminus \mathcal{S}_{L}^+)\leq  e^{-\frac{1}{2} (\log L)^2}.
\end{equation}
\end{lemma}
\begin{proof}
 Let $\mathbf{P}$ be the law of the nearest-neighbor symmetric simple random walk on $\mathbb{Z}$, and
 $(S_i)_{i\in \mathbb{N}}$ be its trajectory with $S_0=0$.
 Since any trajectory of this simple random random walk has the same mass, we have
 \begin{align}
 U_L(\mathcal{S}_{L}\setminus \mathcal{S}_{L}^+)&=\mathbf{P}\Big[\exists i \in \lint 0, L \rint: S_i+m\leq 0 \vert S_L=0\Big] \nonumber\\
 &\leq L^{\frac{1}{2}}\mathbf{P}\Big[\min_{i \in \lint 0, L \rint} S_i \leq -m, S_L=0\Big]  \nonumber\\
 &\leq 2L^{\frac{1}{2}}\mathbf{P}[S_L\leq -m]  \nonumber\\
 &\leq  e^{-\frac{1}{2} (\log L)^2}, \numberthis 
 \end{align}
 which vanishes as $L$ tends to infinity. 
 The first inequality uses $\mathbf{P}[S_L=0]\geq  L^{-1/2}$,  for all $L$ sufficiently large. The second inequality uses $$\mathbf{P}\Big[\min_{i \in \lint 0, L \rint} S_i \leq -m, S_L=0\Big] \leq 2\mathbf{P}\Big[ S_L \leq -m\Big].$$ 
 In the last inequality,  we use the inequality, $\sqrt{2\pi} n^{n+\frac{1}{2}}e^{-n} \leq n! \leq e n^{n+\frac{1}{2}}e^{-n}$ for all  $ n\geq 1$, to obtain
 \begin{align*}
 \mathbf{P}[S_L \leq -m]&\leq (L-m+1)    {L\choose \frac{L+m}{2}}2^{-L}\leq  
 (L-m+1) e^{-(\log L)^2}.
 \end{align*}
 \end{proof}

\subsection{Compare the polymer pinning dynamics to simple exclusion process.} 
The graphical construction mentioned at the beginning of Appendix \ref{proof of the lemma:highestpoint} allows to
couple the three dynamics $(\sigma_t^{\bwedge, 0})_{t\geq 0}$, $(\eta_t^{\bwedge})_{t\geq 0}$ and $(\eta_t^{U_L})_{t\geq 0}$ such that for  all $t\geq 0$,
 \begin{equation}
 \sigma_t^{\bwedge, 0}\geq \eta_t^{\bwedge}\geq \eta_t^{U_L}. \label{monotonicity of the dynamics of maximal path and equilibrium path}
 \end{equation}
Let $ P_t^{\bwedge, -}(\cdot)\colonequals\mathbb{P}(\eta_t^{\bwedge}=\cdot)$ and $P_t^{\bwedge,0}(\cdot)\colonequals\mathbb{P}(\sigma_t^{\bwedge, 0}=\cdot)$. 
Intuitively, the distribution of $\sigma_t^{\bwedge,0}$ is close to that of $\eta^{\bwedge}_t$ for all $t\geq0$.

\begin{lemma} \label{the difference between two maximal path following different parameter}
For any given $\epsilon>0$ and all $L$ sufficiently large, we have
\begin{equation}
\sup_{0\leq t \leq t_{\delta} } \Vert P_t^{\bwedge, 0} -P_t^{\bwedge, -}\Vert_{\TV}\leq \epsilon.
\end{equation}
\end{lemma}

\begin{proof}
 By \cite[Theorem 2.4]{lacoin2016mixing}, for any given $\epsilon>0$ and $t\geq t_{\delta/2}$, if $L$ is sufficiently large, we have 
 \begin{equation}
 \Vert P_t^{\bwedge, -} -U_L \Vert_{\TV}\leq\epsilon.  \label{the distance between the dynamics starting with the maximal path and the equilibrium respectively}
 \end{equation}
 Moreover, by Proposition \ref{characterization of total variation distance} and  the monotonicity in (\ref{monotonicity of the dynamics of maximal path and equilibrium path}),  we obtain
 \begin{align*}
 \sup_{0\leq t \leq t_{\delta} } \Vert P_t^{\bwedge, 0} -P_t^{\bwedge, -}\Vert_{\TV} &\leq \mathbb{P}\Big[ \exists t\in [0, t_{\delta}]: \sigma_t^{\barwedge, 0}\neq \eta_t^{\barwedge}\Big]\\
 &\leq \mathbb{P}\Big[\exists t \in [0, t_{\delta}]: \min_{x\in \lint 0, L \rint}\eta_t^{\bwedge}(x)\leq 0 \Big]\\
 &\leq \mathbb{P}\Big[ \exists t \in [0, t_{\delta}]:\min_{x\in \lint 0, L \rint} \eta_t^{U_L}(x)\leq 0\Big].
  \numberthis  \label{the difference between two maximal path following different parameter_in the proof}
 \end{align*}
 The second inequality is based on the fact that in the coupling if $\sigma_t^{\bwedge, 0}\neq \eta_t^{\bwedge}$, there must exist $x\in \lint 0, L\rint$ satisfying $\eta_s^{\bwedge}(x)=0$ for some $s\in [0, t]$. The third inequality uses the monotonicity of the dynamics, \textit{i.e.} $\eta^{\bwedge}_t\geq \eta^{U_L}_t$ for all $t\geq 0$. The last term in (\ref{the difference between two maximal path following different parameter_in the proof}) vanishes as $L$ tends to infinity, which follows exactly as that in (\ref{upper bound using occupation time}) of Lemma \ref{flippable corners for the equilibrium}, using occupation time (\ref{occupation time}), strong Markov property and Lemma \ref{upper bound for the lifted up system to hit the zero level line}. 
 \end{proof}

 Since $P^{\bwedge,-}_t$ is close to $U_L$ for all $t\geq t_{\delta/2}$, we can use the information of $U_L$ to give an upper bound for the highest point of $\sigma_t^{\bwedge,0}$. 

\begin{proof}[Proof of Lemma \ref{Lemma for controlling the highest point}]  
By triangle inequality, Lemma \ref{the difference between two maximal path following different parameter} and (\ref{the distance between the dynamics starting with the maximal path and the equilibrium respectively}),  for $t\in [t_{\delta/2}, t_{\delta}]$, if $L$ is sufficiently large, we have
\begin{equation}
\Vert P_t^{\bwedge, 0}-U_L\Vert_{\TV} \leq 2\epsilon. \label{total variation distance for the dynamics never hits the x-axis}
\end{equation} 
By (\ref{total variation distance for the dynamics never hits the x-axis}), for every $t \in [t_{\delta/2}, t_{\delta}]$ and $L$ sufficiently large,  we obtain
 \begin{equation*}
 \mathbb{P}\Big[\barh(t)\geq 2L^{\frac{1}{2}}(\log L)^2\Big]  \leq U_L \bigg(\sup_{x\in \lint 0, L \rint }\zeta_x\geq 2L^{\frac{1}{2}}(\log L)^2 , \zeta\in \mathcal{S}_L\bigg) +      \Vert P_t^{\bwedge, 0}-U_L\Vert_{\TV}\leq 3\epsilon,
 \end{equation*}
where the first term in the right hand side vanishes as $L$ tends to infinity, whose proof is the same as Lemma \ref{upper bound for the lifted up system to hit the zero level line}. Since $\epsilon>0$ is arbitrary, we finish the proof.
 \end{proof} 

\section{Spin system.}\label{spinsystem}
To deduce Proposition \ref{Peres-Winkler inequality} from \cite[Theorem 1.1]{peres2013can},  we construct a monotone system  $\langle  \Omega_L^*, S, V_L, \mu_L^* \rangle$ which is the same as the Glauber dynamics of the polymer pinning model.

\begin{figure}[h]
\begin{tikzpicture}[scale=0.7]
\draw[step=1cm,lightgray, ultra thin] (0,0) grid (12, 7);
  \draw (0,0) -- (6,6);
   \draw (6,6) -- (12,0);
 \filldraw[fill=blue!40!white, draw=blue, dashed]  (6,6) -- (5,7)--(0,2)--(1,1)--(6,6); 
       \filldraw[fill=blue!40!white, draw=blue,dashed] (6,6)--(7,7)--(12, 2)--(11, 1)--(6,6);      
        \node [blue] at (5, 6) {$\ominus$};\node [blue] at (4, 5) {$\ominus$}; \node [blue] at (3, 4) {$\ominus$};\node [blue] at (2, 3) {$\ominus$}; \node [blue] at (1, 2) {$\ominus$};
        \node [blue] at (7, 6) {$\ominus$};     \node [blue] at (8, 5) {$\ominus$};      \node [blue] at (9, 4) {$\ominus$};      \node [blue] at (10, 3) {$\ominus$};      \node [blue] at (11, 2) {$\ominus$};    
        \draw[blue, dashed](1,3)--(2,2); \draw[blue, dashed](3,3)--(2,4); \draw[blue, dashed](3,5)--(4,4); \draw[blue, dashed](5,5)--(4,6);\draw[blue, dashed](7,5)--(8,6); \draw[blue, dashed](8,4)--(9,5);\draw[blue, dashed](9,3)--(10,4); \draw[blue, dashed](10,2)--(11,3);
     
   \draw[black,thick, dashed] (2,0)--(7,5);
   \draw[black,thick, dashed] (4,0)--(8,4);
   \draw[black,thick, dashed] (6,0)--(9,3);
   \draw[black,thick,dashed] (8,0)--(10,2);
   \draw[black,thick,dashed] (10,0)--(11,1);
   \draw[black,thick,dashed] (1, 1)--(2,0);
   \draw[black,thick,dashed] (2,2)--(4,0);
   \draw[black,thick,dashed] (3,3)--(6,0);
    \draw[black,thick,dashed] (4,4)--(8,0);
     \draw[black,thick,dashed] (5,5)--(10,0);      
     \draw[blue,thick,dashed] (7,7)--(12,2);
     \draw[black,thick,dashed] (0,0)--(6,6);
     \draw[black,thick,dashed] (6,6)--(12,0);

      \filldraw[fill=green!40!white, draw=black, dashed] (1,1) -- (2,2) -- (3,1) -- (2,0) -- (1,1); 
      \filldraw[fill=green!40!white, draw=black,dashed] (3,1) -- (4,0) -- (5,1) -- (4,2) -- (3,1); 
       \filldraw[fill=green!40!white, draw=black,dashed] (5,1) -- (6,0) -- (7,1) -- (6,2) -- (5,1); 
          \filldraw[fill=green!40!white, draw=black,dashed] (7,1) -- (8,0) -- (9,1) -- (8,2) -- (7,1); 
             \filldraw[fill=green!40!white, draw=black,dashed] (9,1) -- (10,0) -- (11,1) -- (10,2) -- (9,1);
        \filldraw[fill=red!40!white, draw=black,dashed] (0,0) -- (1,1) -- (2,0)-- (0,0);
              \filldraw[fill=red!40!white, draw=black,dashed] (2,0) -- (3,1) -- (4,0)-- (2,0);
              \filldraw[fill=red!40!white, draw=black,dashed] (4,0) -- (5,1) -- (6,0)-- (4,0);
             \filldraw[fill=red!40!white, draw=black,dashed] (6,0) -- (7,1) -- (8,0)-- (6,0);
             \filldraw[fill=red!40!white, draw=black,dashed] (8,0) -- (9,1) -- (10,0)-- (8,0);
              \filldraw[fill=red!40!white, draw=black,dashed] (10,0) -- (11,1) -- (12,0)-- (10,0);
     \draw[blue , line width=0.6mm, -] (0,0) -- (3,3);
            \draw[blue , line width=0.6mm, -] (4,2) -- (3,3);
            \draw[blue , line width=0.6mm, -] (4,2) -- (5,3);
            \draw[blue , line width=0.6mm, -] (7,1) -- (5,3);
            \draw[blue , line width=0.6mm, -] (7,1) -- (8,0);
             \draw[blue , line width=0.6mm, -] (8,0) -- (9,1);
              \draw[blue , line width=0.6mm, -] (9,1) -- (10,2)--(12,0);          
              \draw[blue , line width=0.5mm, dashed] (7,1) -- (8,2)--(9, 1);                     
       \node [red] at (2,1) {$ \mathbf{\oplus} $};
        \node [red] at (4,1) {$\oplus$};
         \node [red] at (6,1) {$\oplus$};
          \node [blue] at (8,1) {$\ominus$};
           \node [red] at (10,1) {$\oplus$};
            \node [red] at (3.1,2.1) {$\oplus$};
            \node [red] at (5.1,2.1) {$\oplus$};
          \node [red] at (1,0.4) {$\oplus$};   \node [red] at (3,0.4) {$\oplus$};   \node [red] at (5,0.4) {$\oplus$};   \node [red] at (7,0.4) {$\oplus$};   \node [red] at (9,0.4) {$\oplus$};   \node [red] at (11,0.4) {$\oplus$};   
          \node [blue] at (4.1,3.1) {$\ominus$};
             \node[blue] at (5.1,4.1) {$\ominus$};
                \node [blue] at (6.1,5.1) {$\ominus$};
                   \node [blue] at (6.1,3.1) {$\ominus$};
                      \node [blue] at (7.1,2.1) {$\ominus$};
                         \node [blue] at (7.1,4.1) {$\ominus$};
                            \node [blue] at (8.1,3.1) {$\ominus$};
                               \node [blue] at (9.1,2.1) {$\ominus$};
                               
         \node  at (0, -0.3) {0};                      
        \node  at (2, -0.3) {2}; \node  at ( -0.15, 2) {2};
             \node  at (4, -0.3) {4};  \node  at (-0.15,  4) {4};
                  \node  at (6, -0.3) {6};    \node  at (-0.15, 6)  {6};
                       \node  at (8, -0.3) {8};    
                        \node  at (10, -0.3) {10};
                             \node  at (12, -0.3) {12};
                       \node at(5.8, 2.45){$\xi$};
 \draw[thick,->] (0,0) -- (12.5,0) node[anchor=north west]{x};
\draw[thick,->] (0,0) -- (0,7.3) node[anchor=south east] {y};                      
\end{tikzpicture}
\caption { \label{spinsystemfig}An example shows the equivalence between the polymer pinning model and the spin system with $L=12$.
The blue path $\xi$  is an element of $\Omega_L$. This configuration in the spin system is denoted as $\xi^*$, and its probability measure is $\mu(\xi)$.
  The corner at $x=8$ of thick blue path $\xi$  flips with rate $1/(1+\lambda)$ to the dashed blue corner, while the spin $\ominus$ at the green square centered at $(8, 1)$ flips to $\oplus$ with rate $1/(1+\lambda)$.}
\end{figure}
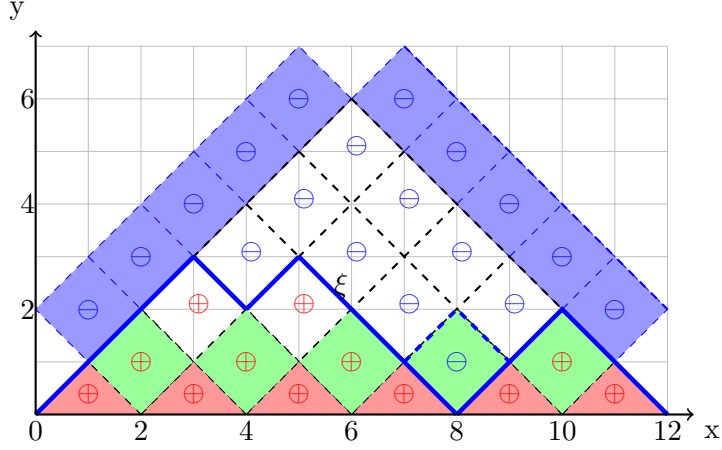

For $(x, z)\in \mathbb{N}^2$, a square with four vertices $\{(x-1,z), (x+1, z), (x, z-1), (x,z+1)\}$ is denoted as $Sq(x,z)$. Recalling $\Theta$ defined in (\ref{the set of spins for censoring}), let $S\colonequals \{\oplus, \ominus\}$ denote the spins, and  $V_L \colonequals \{ Sq(x,z): \forall (x, z)\in \Theta \}$ denote the set of all sites, which consists of all green or white squares shown Figure \ref{spinsystemfig}.  Each square of $V_L$ is endowed with  $\oplus$ or $\ominus$. Moreover, we give a natural order for the spins, say, $\ominus \leq \oplus$.  For any given $\xi \in \Omega_L$, every square $Sq(x,z)$ lying under the path $\xi$ is endowed with $\oplus$, while every square $Sq(x, z)$ lying above $\xi$ is endowed $\ominus$. 
 This spin configuration is denoted as $\xi^*$.  For $\xi, \xi'\in \Omega_L$, 
$\xi\leq \xi'$ if and only if $\xi^*\leq \xi'^{*}$.
 Let $\Omega_L^*\colonequals \{\xi^*, \xi\in \Omega_L \}$ and $\mu_L^*(\xi^*) \colonequals \mu(\xi)$.

 For convenience of describing the Glauber dynamics of spin system, we introduce two fixed boundary conditions. We  assign a negative spin $\ominus$ to each square $Sq(x, z)$ where $$\Big\{(x, z): x\in \lint 1, L/2-1\rint \cup \lint L/2+1, L-1 \rint , z=L/2+1-\vert x-L/2\vert \Big\}.$$ These are the blue squares shown in Figure \ref{spinsystemfig}.  In addition, we also introduce a positive boundary condition.  A triangle with three vertices $\{(x-1, z), (x+1, z), (x, z+1)\}$ is denoted as $Tr(x, z)$ for $(x,z)\in \mathbb{N}^2$. We assign a positive spin $\oplus$ to each triangle $Tr(x,0)$ for all $x\in \lint 1, L-1\rint \setminus 2\mathbb{N}$. These are the red triangles shown in Figure \ref{spinsystemfig}.
 We say that two spins are neighbors if the squares or triangles, they lie, share an edge.  We use the same exponential clocks and uniform coins $\Tau^{\uparrow}$, $\Tau^{\downarrow}$, $\mathcal{U}^{\uparrow}$, and $\mathcal{U}^{\downarrow}$ define in Subsection \ref{graphical construction} to describe the dynamics of the spin system.

Given $\Tau^{\uparrow}$, $\Tau^{\downarrow}$, $\mathcal{U}^{\uparrow}$ and $\mathcal{U}^{\downarrow}$,  we construct, in a deterministic way,  $(\sigma_t^{\xi^*})_{t \geq 0}$  the Glauber dynamics of the spin system  starting with $\xi^* $ with parameter $\lambda$. The trajectory $(\sigma_t^{\xi^*})_{t\geq 0}$ is c\`{a}dl\`{a}g with $\sigma_0^{\xi^*}=\xi^*$ and is constant in the intervals, where the clock processes are silent. 

When the clock process $\Tau^{\uparrow}_{(x, z)}$ rings at time $t=\Tau_{(x, z)}^{\uparrow}(n)$  for $n\geq 1$, we update the configuration $\sigma_{t^-}^{\xi^*}$ as follows:
\begin{itemize}
 \item if the spin in the square $Sq(x, z)$ is $\ominus$, and has two neighbors with $\oplus$ spins,  and $z=1$, and $\mathcal{U}^{\uparrow}_{(x,z)}(n)\leq \frac{1}{1+\lambda}$, we let the spin in the square $Sq(x, z)$ change to $\oplus$ at time $t$, and the other spins remain unchanged;
\item if the spin in the square $Sq(x, z)$ is $\ominus$, and has two neighbors with $\oplus$ spins, and $z> 1$, and $\mathcal{U}^{\uparrow}_{(x,z)}(n)\leq 1/2$, we let the spin in the square $Sq(x, z)$ change to $\oplus$.
\end{itemize}
 If these two conditions aforementioned are not satisfied, we do nothing.

When the clock process $\Tau^{\downarrow}_{(x, z)}$ rings at time $t=\Tau_{(x,z)}^{\downarrow}(n)$ for $n\geq 1$, we update the configuration $\sigma_{t^-}^{\xi^*}$ as follows:
\begin{itemize}
\item  if the spin in the square $Sq(x, z)$ is $\oplus$, and has two neighbors with $\ominus$ spins,  and $z=1$, and $\mathcal{U}^{\downarrow}_{(x,z)}(n)\leq \frac{\lambda}{1+\lambda}$, we let the spin in the square $Sq(x, z)$ change to $\ominus$ at time $t$, and the other spins remain unchanged;
\item  if the spin in the square $Sq(x, z)$ is $\oplus$, and has two neighbors with $\ominus$ spins,  and $z\geq 2$, and $\mathcal{U}^{\downarrow}_{(x,z)}(n)\leq 1/2$, we let the spin in the square $Sq(x, z)$ change to $\ominus$ at time $t$, and the other spins remain unchanged.
\end{itemize}
 If these two conditions aforementioned are not satisfied, we do nothing.

We can see that   $\langle  \Omega^*, S, V, \mu^* \rangle$ is a monotone system in the sense of \cite[Section 1.1]{peres2013can},
whose Glauber dynamics is the same as that of the polymer pinning model.

\end{appendix}
\bibliographystyle{alpha}
\bibliography{/home/shangjie/Desktop/probability-paper-projects/library}

\end{document}